\documentclass{amsart}
\usepackage{psfrag, comment, xypic}
\usepackage{amsmath,amsfonts,amssymb,amsthm,epsfig}

\newcommand{\R}{{\mathbb R}}
\newcommand{\Z}{{\mathbb Z}}
\newcommand{\Q}{{\mathbb Q}}

\newcommand{\CH}{{\mathcal H}}

\def\cL{\mathcal{L}}
\def\ts{\hspace{-.02cm}}
\def\qe{$\diamond$}
\def\dontshow#1
{ }
\newcommand{\tubelattice}{\mathcal{L}}
\newcommand{\Cmn}{\mathcal{C}}
\newcommand{\semifield}{\mathbb F}
\newcommand{\bW}{\overline{W}}

\newtheorem{Theorem}{Theorem} 
\newtheorem{Proposition}[Theorem]{Proposition} 
\newtheorem{Conjecture}[Theorem]{Conjecture} 
\newtheorem{Lemma}[Theorem]{Lemma}

\newtheorem{Corollary}[Theorem]{Corollary}

\theoremstyle{definition} 
\newtheorem{Example}[Theorem]{Example}
\newtheorem{Remark}[Theorem]{Remark}

\begin{document}

\title{A periodicity theorem for the octahedron recurrence}
\thanks{$\dagger$ An earlier version of this work has circulated under the name ``A coboundary category defined using the octahedron recurrence''}
\author{Andr\'e Henriques}

\address{	
Andr\'e Henriques\\
Mathematisches Institut \\
Westf\" alische Wilhelms-Universit\" at \\
Einsteinstr. 62\\
48149\\ M\" unster\\
Germany}
\email{andrhenr@math.uni-muenster.de}

%\date{\today}
%\begin{abstract} \end{abstract}
%\tableofcontents

\maketitle

\begin{abstract}
The octahedron recurrence lives on a 3-dimensional lattice and is given by
\(
f(x,y,t+1)=\big(f(x+1,y,t)f(x-1,y,t)+f(x,y+1,t)f(x,y-1,t)\big)\big/f(x,y,t-1)
\)
In this paper, we investigate a variant of this recurrence which lives in a
lattice contained in $[0,m] \times [0,n] \times \mathbb R$.
Following Speyer, we give an explicit non-recursive formula for the values of this
recurrence and use it to prove that it is periodic of period $n+m$.
We then proceed to show various other hidden symmetries satisfied this bounded octahedron recurrence.
\end{abstract}

\section{Introduction}

In this paper we investigate a variant of the octahedron recurrence of Robbins-Rumsey \cite{RR86} called the {\em bounded octahedron recurrence}. 
It was first described by Kamnitzer and the author in \cite{HK05}, 
where it was used to relate the commutativity isomorphism for $\mathfrak{gl}(n)$-crystals with the Sch\" utzenberger involution on Young tableaux.

The bounded octahedron recurrence takes place on the lattice\dontshow{Lde}
\begin{equation}\label{Lde}
\cL:=\big\{(x,y,t)\in\Z^3\,\big|\,
0\le x\le m,\, 
0\le y\le n,\,
x+y+t\equiv 0 \mod 2
\big\}
\end{equation}
and is best described by the following figures:\dontshow{iaf}

{   \psfrag{max(a+c,b+d)-e}{$(ac + bd)/e$}
   \psfrag{a}{$a$}
   \psfrag{b}{$b$}
   \psfrag{c}{$c$}
   \psfrag{d}{$d$}
   \psfrag{e}{$e$}\psfrag{a+c-e}{$ac/e$}\psfrag{a+b-e}{$ab/e$}
\begin{equation}\label{iaf}
\begin{matrix}\epsfig{file=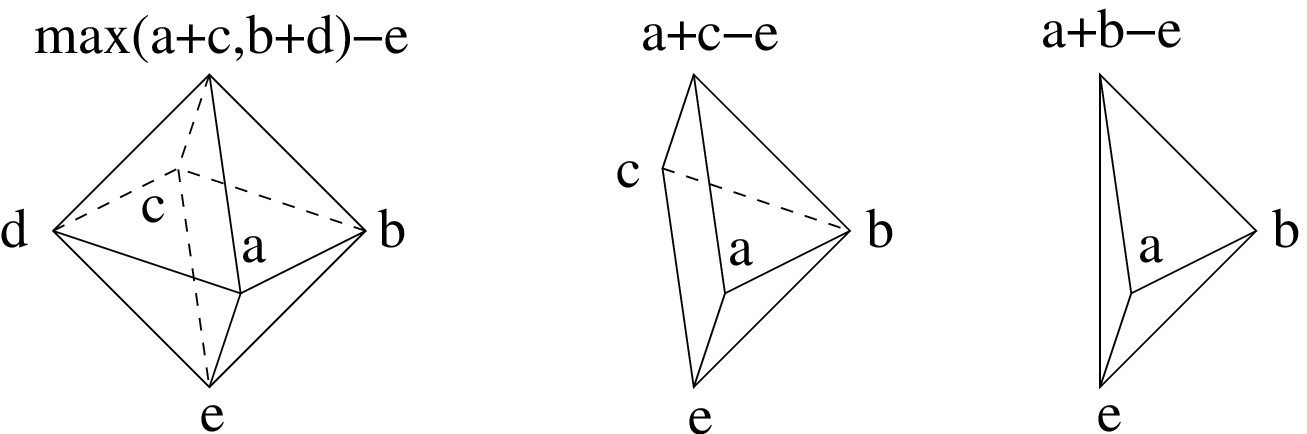,height=2cm}\end{matrix}
\end{equation}
%\centerline{The bounded octahedron recurrence.}
}

\noindent
We think of the first two coordinates as space, and of the third one as time.
The values at the points of $\cL$ with higher third coordinate (points in the future) are then computed from the values at the points of $\cL$ 
with lower third coordinate (points in the past).
The feature that distinguishes this recurrence from \cite{RR86} is that the space coordinates are now bounded.
One then has two additional rules that describe the recurrence for points on the boundary, and for points on the corners respectively.

The original (unbounded) octahedron recurrence was studied by Fomin and Zelevinsky \cite{FZ02} as an example of the `Laurent phenomenon'.
Generalizing \cite{RR86}\footnote{Fomin and Zelevinsky attribute this to \cite{MRR83}, but we believe that it is best attributed to \cite{RR86}.}, 
they proved that the value at a future point is always a Laurent polynomial in terms of the initial data.
Speyer \cite{Spe04} then refined their results by showing that the monomials in the above Laurent polynomials are in bijection
with the set of perfect matchings of certain graphs, and that the coefficients are all 1.

Following Speyer, we give a non-recursive formula for the bounded octahedron recurrence using perfect matchings.
Namely, we prove that the solution of the recurrence at some future point is given as a sum over matchings of a certain graph.  
Our formula is very similar to that of Speyer. 
However, there is one important distinction: 
in the unbounded case, the size of the graph grows linearly as you move your point into the future.  
In the bounded case, this graph grows linearly for a while but then shrinks again.  
In particular, we use this to prove that the octahedron recurrence is periodic of period $n+m$:\dontshow{pPt}
\begin{Theorem}\label{pPt} 
Let $f$ be a function satisfying the bounded octahedron recurrence, then there is a constant $c$ such that for any point $(x,y,t)\in \cL$ 
we have the relation
\begin{equation*}
f(x,y,t)=c\cdot f(m-x,n-y,t-m-n).
\end{equation*}
\end{Theorem}
This result is very reminiscent of Fomin and Zelevinsky's theorem about the periodicity of $Y$-systems \cite{FZ03}.

As an application of our main theorem, we deduce certain remarkable identities between the bounded octahedron recurrence in various domains.
Finally, we conjecture a similar periodicity for the bounded cube recurrence, which is a similar variant of the recurrence studied in \cite{Pro01}, \cite{FZ02}, \cite{CS04}.

\subsection{Acknowledgements}
I would like to thank Alexandre Goncharov, Allen Knutson, and Dylan Thurston for helpful conversations.
I reserve special thanks for Joel Kamnitzer and David Speyer, who have greatly influenced me as I was doing this work.
I also thank David Speyer for a careful reading of this document.

\section{The Bounded Octahedron Recurrence}\label{s:bo}

Throughout\dontshow{s:bo} this paper, we will be studying functions with values in a semifield $ \semifield $.  
A \textit{semifield} is a set $ \semifield $ along with two operations called addition and multiplication, such that:
\begin{enumerate}
\item addition is commutative and associative,
\item multiplication makes $ \semifield $ into an abelian group\footnote{Note that fields are not semifields since they contain a zero element.},
\item multiplication distributes over addition.
\end{enumerate}

There are two main classes of examples of semifields.  The first are positive parts of ordered fields such as $ \R_{>0} = \{ x \in \R : x > 0 \} $ or 
$ \Q_{>0} $ under the usual operations.  The second are the tropical semifields $ \Z_t, \Q_t, \R_t $ where addition is $ \max $ and multiplication 
is $ + $.
 
Fix $m,n\in \Z_{>0}$.
Let us call {\em space-time} the space $Y=[0,m]\times[0,n]\times \R$. The first two coordinates represent ``space'' and the last one is ``time''. 
In $Y$, we have the lattice 
$\tubelattice=\{(x,y,t)\in\Z^3\cap Y :\, x+y+t \text{ is even}\}$ on which the recurrence will take place.  
It is the set of vertices of a tiling of $Y$ by tetrahedra, octahedra,
$1/2$-octahedra, and $1/4$-octahedra as shown in (\ref{tst}). 
The tetrahedra are given by
\begin{align*}
conv\{&(x,y,t),(x+1,y+1,t),(x+1,y,t+1),(x,y+1,t+1)\},\quad  x+y+t \text{ even,}\\
conv\{&(x+1,y,t),(x,y+1,t),(x,y,t+1),(x+1,y+1,t+1)\},\quad  x+y+t \text{ odd,}
\end{align*}
%\intertext{
while the octahedra, 1/2-octahedra and 1/4-octahedra are given by
\begin{equation*}
\begin{split}
Y\cap conv\{(x\ts+\ts 1,y,t),(x,y\ts+\ts 1,t),(x,y,t\ts+\ts 1),(x\ts-\ts 1,y,t),(x,y\ts-\ts 1,&t),(x,y,t\ts-\ts 1)\},\\ 
&x+y+t \text{ odd.}
\end{split}
\end{equation*}

\begin{gather}\label{tst}
\begin{matrix}\epsfig{file=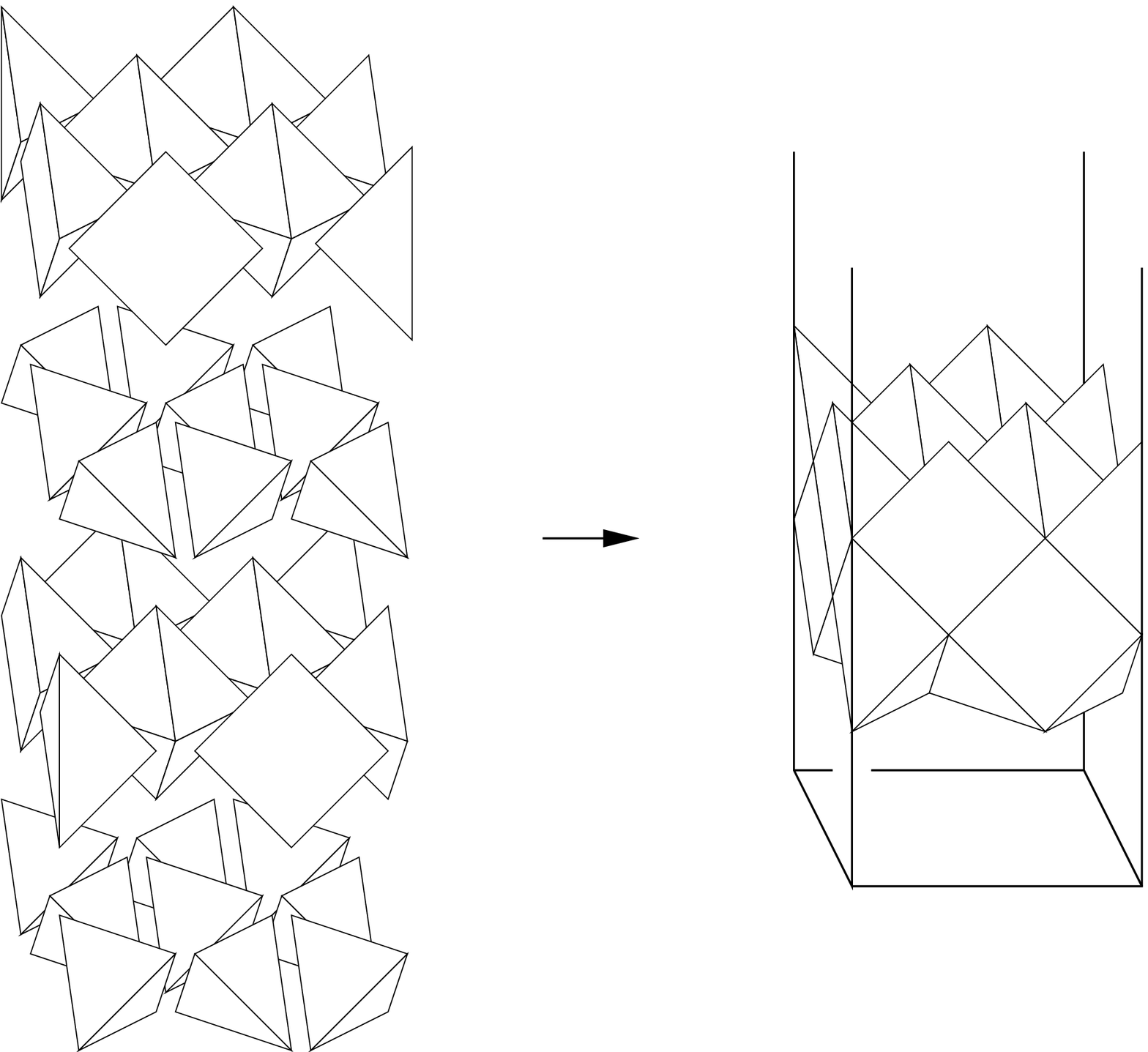,height=6cm}\end{matrix}\\
\nonumber\text{The tiling of space-time.}
\end{gather}
\vspace{.3cm}

\dontshow{tst}
A {\em section} is a subcomplex $S$ of the 2-skeleton 
of the above tiling which contains exactly one point over each $ (x,y) $.  
In particular, $ S $ is the graph $ S = \{ (x,y,h(x,y)) \} $ of a continuous 
map $h:[0,m]\times[0,n]\to\R$. 
A point $ (x,y,t) \in \tubelattice $ is said to be in the {\em future}
of a section $ S $ if there exists $ (x, y, t') \in S $ with $ t' \le t $.

A {\em state} of a subset $ A \subset Y $ is an $\semifield$-valued function 
$f:A\cap\tubelattice\to \semifield$. In particular we may speak of the state of a section.
The state $f$ of a section $S$ determines 
the state (again denoted by $f$) of the set of all points in its future,
according to the {\em bounded octahedron recurrence}:\dontshow{octrec}

\begin{align}
\nonumber f(x,y,t\ts+\ts 1)=
&\big(f(x\ts+\ts 1,y,t)f(x\ts-\ts 1,y,t) + f(x,y\ts+\ts 1,t)f(x,y\ts-\ts 1,t)\big)\big/f(x,y,t\ts-\ts 1) \\
\nonumber &\phantom{f(x\ts+\ts 1,y,t)f(x\ts-\ts 1,y,t)/f(x,y,t\ts-\ts 1) \hspace{.8cm}} \text{if } 0\ts<\ts x\ts<\ts m,\, 0\ts<\ts y\ts<\ts n,\\
\nonumber &f(x\ts+\ts 1,y,t)f(x\ts-\ts 1,y,t)/f(x,y,t\ts-\ts 1)          \hspace{.8cm} \text{if } 0\ts<\ts x\ts<\ts m, y=0\text{ or }n,\\
\nonumber &f(x,y\ts+\ts 1,t)f(x,y\ts-\ts 1,t)/f(x,y,t\ts-\ts 1)          \hspace{.8cm} \text{if } 0\ts<\ts y\ts<\ts n, x=0\text{ or }m,\\
\label{octrec} &f(x\ts+\ts 1,y,t)f(x,y\ts+\ts 1,t)/f(x,y,t\ts-\ts 1)     \hspace{.8cm} \text{if } (x,y)=(0,0),\\
\nonumber &f(x\ts+\ts 1,y,t)f(x,y\ts-\ts 1,t)/f(x,y,t\ts-\ts 1)          \hspace{.8cm} \text{if } (x,y)=(0,n) ,\\
\nonumber &f(x\ts-\ts 1,y,t)f(x,y\ts+\ts 1,t)/f(x,y,t\ts-\ts 1)          \hspace{.8cm} \text{if } (x,y)=(m,0),\\
\nonumber &f(x\ts-\ts 1,y,t)f(x,y\ts-\ts 1,t)/f(x,y,t\ts-\ts 1)          \hspace{.8cm} \text{if } (x,y)=(m,n). 
\end{align}
\vspace{.2cm}

\noindent So we have one rule if our point is in the interior (this is the original octahedron recurrence in \cite{RR86}), 
another rule if it lies on a wall, and a third 
if it lies on a vertical edge.  These rules can be seen in (\ref{iaf}).

Note that the above formulas are invertible, and are equal to their own inverse. Indeed, $e'=(ac +bd)/e$ if and only if
$e=(ac + bd)/e'$, and similarly for the other ones. 
The state of a section therefore also determines the state of all the
points in its past.
If $f:S\cap\cL\to \semifield$ is a state, we will often abuse notation, and
denote by $f$ its extension to the whole of $\cL$.

\section{Speyer's formula}\subsection{The unbounded case}\label{3p1}

In this section, we recall the main result of Speyer \cite{Spe04}. 
We shall assume for a moment that $Y=\R^2\times \R$ instead of $[0,m]\times[0,n]\times\R$, and that
$\tubelattice$ is the whole lattice $\{(x,y,t)\in\Z^3\,:\, x+y+t \text{ is even}\}$.

Given a section $S$, along with a state $f$ and a point $(x_0,y_0,t_0)$ in its future, the goal is to provide 
an explicit non-recursive formula for $f(x_0,y_0,t_0)$ in terms of $f|_S$.
Since the state of a point only influences the state of its neighbors, 
$f(x_0,y_0,t_0)$ is entirely determined by the restriction of $ f $ to the intersection of $ S $ with the {\em light cone}: 
\begin{equation*} 
\begin{split}
\Cmn=\Cmn(x_0,y_0,t_0):=\{(x,y,t) : \ &t_0+x_0+y_0\ge t+x+y, \ t_0+x_0-y_0\ge t+x-y,\\
  &t_0-x_0+y_0\ge t-x+y,\ t_0-x_0-y_0\ge t-x-y\}.
\end{split}
\end{equation*} 
%Then, he gives an explicit non-recursive combinatorial formula for $f(x_0,y_0,t_0)$ in terms of the state of 
%$f|_{S\cap\Cmn}$.
Let $W\subset S$ be the closed subcomplex given by\dontshow{gmn}
\begin{equation}\label{gmn}
W=W(S,x_0,y_0,t_0):=\begin{cases} 
\begin{matrix}\overline{S\cap\overset{\circ}{\Cmn}}\end{matrix}\quad&\text{if}\quad(x_0,y_0,t_0)\not\in S \\
\{(x_0,y_0,t_0)\} \quad&\text{if}\quad(x_0,y_0,t_0)\in S,\end{cases}
\end{equation} 
where $\overset{\circ}{\Cmn}$ is the interior of $\Cmn$.
{
\psfrag{S}{$S$}\psfrag{C}{$\mathcal C$}\psfrag{W}{$W$}\psfrag{x}{$(x_0,y_0,t_0)$}
\[
\epsfig{file=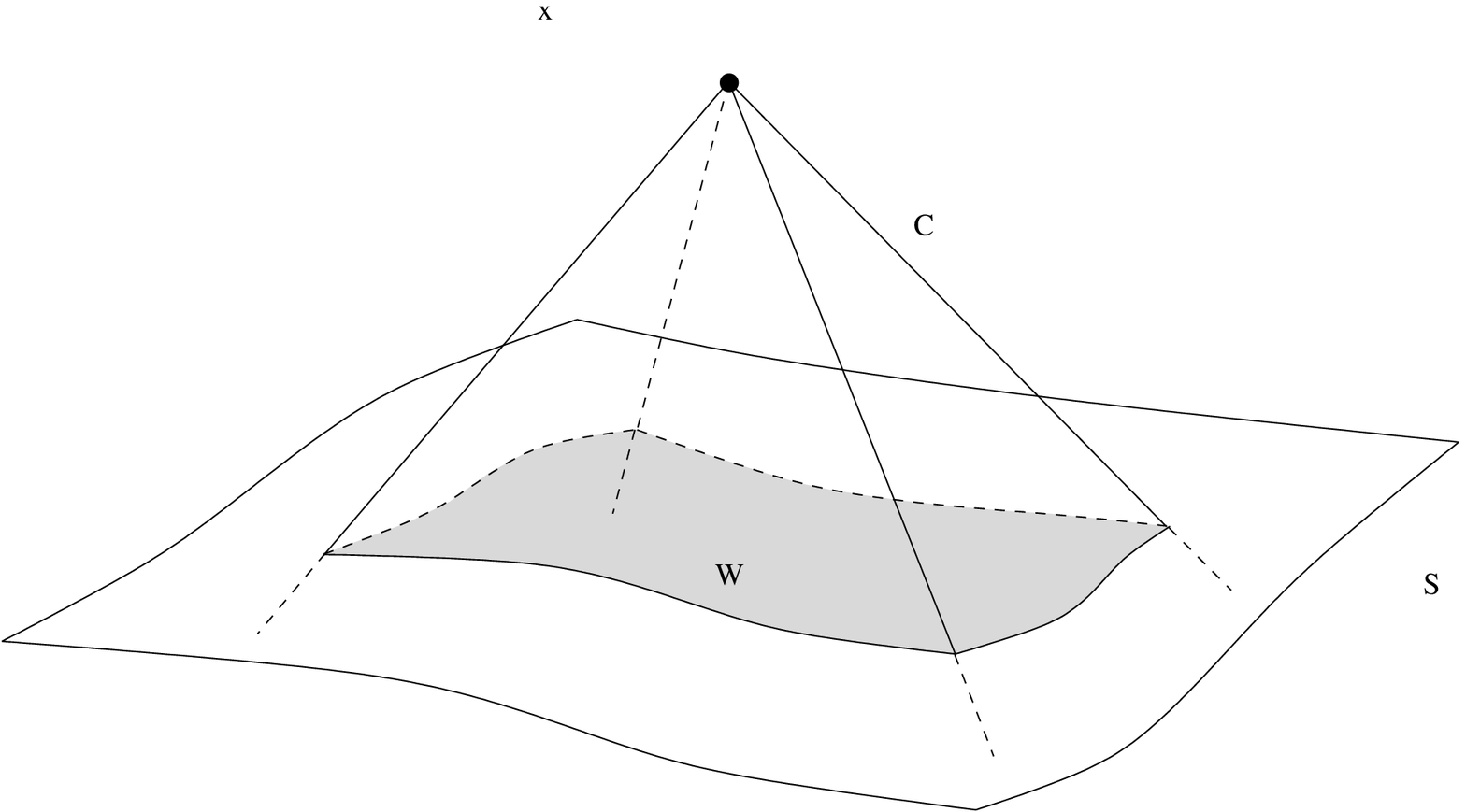,height=2.5cm}
\]
}
We equip $W$ with a cell structure that we now describe. 
Consider the triangulation inherited from the tiling of $Y$. 
An edge $(x,y,t)-(x',y',t')$ is called {\em horizontal} if $t=t'$.
These are the edges whose projections to $\R^2$ are parallel to $x=y$ or $x=-y$.
An edge is called {\em bent} if it's not in $\partial W$, and if the two triangles of $W$ containing it are not coplanar.

The cell structure on $W$ is obtained from the above triangulation by the following two steps.
First delete all the horizontal bent edges (thus creating squares in the corresponding projection onto $\R^2$). 
Then delete the vertices of $\partial W$ that belonged to these horizontal bent edges. 
Note that this creates new boundary edges that are not straight.

\begin{Example}\label{extet}
We illustrate the complex $W$ corresponding to a particular section $S$ and point $(x_0,y_0,t_0)$.
The section is represented first in 3 dimensions and then via its projection to $\R^2$.
The shaded part is the intersection with the interior of the cone $\Cmn$. 
{
\psfrag{s}{$S=$}\psfrag{=}{$=$}\psfrag{w}{$W=$}\psfrag{(x,y,t)}{$\scriptstyle (x_0,y_0,t_0)$}
\[
\epsfig{file=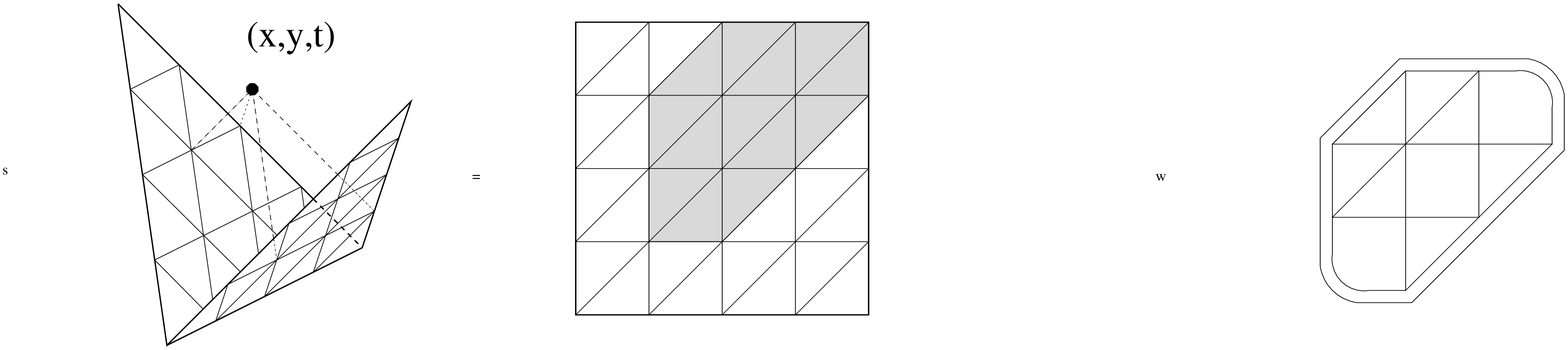,height=2.2cm}
\]
}
\end{Example}

\noindent A {\em matching\footnote{This is not the standard use of the word `matching' in graph theory. 
It is borrowed from Speyer \cite{Spe04} who works with the planar dual of $W$. 
Maybe it would have been more appropriate to call them `dual matchings'.}} $M$ of $W$ is a collection of internal edges of $W$, 
such that every face has exactly one of them on its boundary. 
We will draw the edges $e\in M$ in dotted lines and the edges $e\not\in M$ in solid lines. 
We will also refer to them as dotted edges and solid edges.  
By convention, when we talk about solid edges, we will exclude the boundary edges.

\begin{gather*}
\epsfig{file=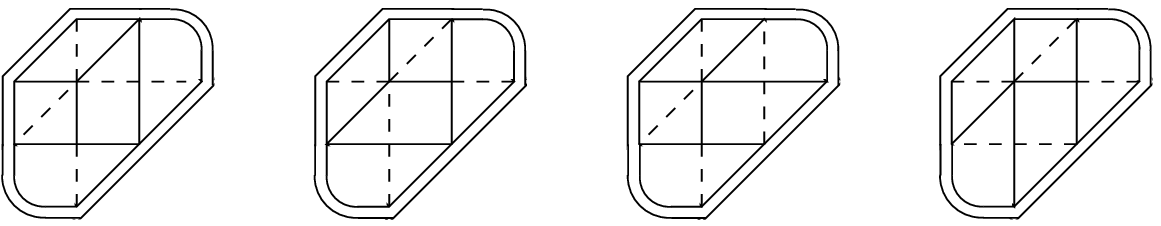,height=1.5cm}\\
\text{The matchings of $W$. }
\end{gather*}

\noindent Suppose we are give a state $f$ of $S$. 
Associated to every matching $M$ of $W$, there is a {\em matching monomial} given by
$$
\mu(M)=\mu(M,f)=
\prod_{(x,y,t)\in W\cap\tubelattice}f(x,y,t)^{k(x,y,t)},$$
where
$$k(x,y,t) =\begin{cases}
\parbox{8cm}{
$1/2\big(\,\#$ of solid edges\ incident to $(x,y,t)$

$\phantom{1/2}-\:\#$ of dotted edges incident to $(x,y,t)\big)$ $-1$}
&
\hspace{-.3cm} \text{if } (x,y,t)\in\overset{\circ}{W},\vspace{.2cm}\\
\parbox{8cm}{
$\big\lceil 1/2\big(\,\#$ of solid edges incident to $(x,y,t)$

$\phantom{1/2}-\: \#$ of dotted edges incident to $(x,y,t)\big)\big\rceil$}
& 
\hspace{-.3cm} \text{if } (x,y,t)\in\partial W,\vspace{.2cm}\\
\hspace{3.5cm}1
&
\hspace{-.3cm} 
\parbox{3cm}{
if $(x,y,t)=$ \hfill 

$\hspace{.32cm} (x_0,y_0,t_0)$.}\end{cases}$$

\noindent
By convention, $\overset{\circ}{W}=\partial W=\emptyset$ in the degenerate case $W=\{(x_0,y_0,t_0)\}$. 
We also recall that the boundary edges don't count as `solid edges'.
Given the above notation, we then have the following result.

\begin{Theorem}
[Speyer \cite{Spe04}]\label{thmoct}
Let $S$ be a section and let $f:S\cap\cL\to\semifield$ be a state (recall that $S$ and $\tubelattice$ are here the unbounded analogs of the notions introduced in Section \ref{s:bo}).
Let $(x_0,y_0,t_0)$ be a point in the future of $S$ such that the corresponding complex $W$ is finite. 
Then the value $f(x_0,y_0,t_0)$ given by the octahedron recurrence is expressed from the initial data by the formula
$$
f(x_0,y_0,t_0)=\sum_{M \in \text{\rm matchings of } W}\mu(M).
$$
\end{Theorem}

\subsection{The bounded case} \label{sec:Tbc}

\dontshow{sec:Tbc} 
We now go back to the situation where $Y=[0,m]\times[0,n]\times \R$, and $\tubelattice\subset Y$ is the corresponding bounded lattice.

In our case of bounded spacetime, we need to consider the following modified light cone which 
``reflects on the walls'': %, illustrated in (\ref{fig:lightcone}):\dontshow{cmn}
\begin{equation}\label{cmn}
\begin{split}
\Cmn=\{(x,y,t)\in Y : \ 
&t_0\ts +\ts x_0\ts +\ts y_0\ge t\ts +\ts x\ts +\ts y,\, \
t_0\ts +\ts x_0\ts -\ts y_0\ge t\ts +\ts x\ts -\ts y,\\
&t_0\ts -\ts x_0\ts +\ts y_0\ge t\ts -\ts x\ts +\ts y,\, \
t_0\ts -\ts x_0\ts -\ts y_0\ge t\ts -\ts x\ts -\ts y,\\
&t\ts +\ts x\ts +\ts y\ge t_0\ts -\ts x_0\ts -\ts y_0,\, \
t\ts +\ts x\ts -\ts y\ge t_0\ts -\ts x_0\ts -\ts (2n\ts -\ts y_0),\\
&t\ts -\ts x\ts +\ts y\ge t_0\ts -\ts (2m\ts -\ts x_0)\ts -\ts y_0, \\
&t\ts -\ts x\ts -\ts y\ge t_0\ts -\ts (2m\ts -\ts x_0)-(2n\ts -\ts y_0)\}.
\end{split}
\end{equation}
In this case the light cone is not a cone at all: it's a rhombic dodecahedron. 
Amazingly, with this new light cone an appropriately modified version of Theorem \ref{thmoct} holds.
{
\psfrag{a}{$(x_0,y_0,t_0)$}
\psfrag{b}{$(m-x_0,n-y_0,t_0-m-n)$}
\begin{gather*}%\label{fig:lightcone}
\begin{matrix}\epsfig{file=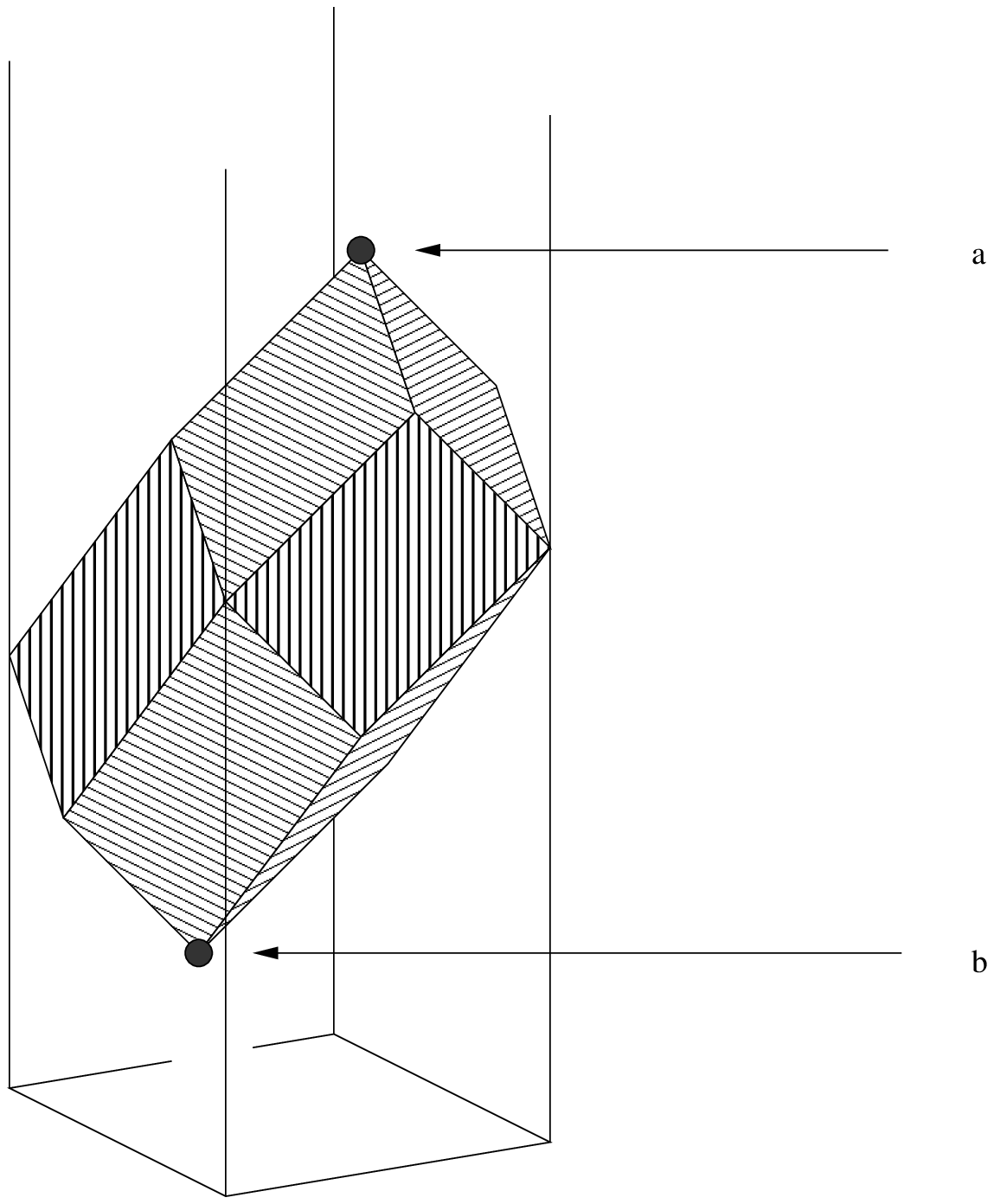,height=5cm}\end{matrix}\hspace{2cm}
\\
\nonumber\text{The modified light cone $\Cmn$.}
\end{gather*}
}

Intuitively, we can think of this in the following manner.  When you want to compute the state of a point $(x_0,y_0,t_0)$ 
\dontshow{fig:lightcone}
which is far away in the future, but not too far, you need to do a big computation since the size of $W$ is large. But 
as your point goes even farther in the future (farther than the diameter of the universe), the complex $W$ starts becoming 
smaller and the computation starts becoming easier. 
When $(x_0,y_0,t_0)$ is exactly $n+m $ steps into the future, then $W$ is reduced to a single point.  
This gives us the periodicity of the recurrence.

We now formulate our version of Theorem \ref{thmoct}. 
Let $S$ be a section, $f$ a state, and $(x_0,y_0,t_0)\in \tubelattice$ a point in the future of $S$. 
We assume $(x_0,y_0,t_0)$ is in the interior of $Y$. 
We will also assume that $\Cmn$ meets $S$ in at least one point.
This is equivalent to having the point $(m-x_0,n-y_0,t_0-m-n)$, which is at the bottom of $\Cmn$, be either in $S$ or in the past of $S$. Let
\begin{equation}\label{subcomplex}
\tilde W=\begin{cases} 
\begin{matrix}\overline{S\cap\overset{\circ}{\Cmn}}\end{matrix}
\quad
&
\parbox{5.5cm}{
if $\hspace{.06cm}(x_0,y_0,t_0)\not\in S\hspace{.06cm}$ and

$\hspace{.38cm} (m\ts -\ts x_0,n\ts -\ts y_0,t_0\ts -\ts m\ts -\ts n)\not\in S$.}\vspace{.1cm}\\
\{(x_0,y_0,t_0)\} \quad&\text{if }\,(x_0,y_0,t_0)\in S,\\
\{(m\ts -\ts x_0,n\ts -\ts y_0,t_0\ts -\ts m\ts -\ts n)\} \quad&\text{if }\,(m\ts -\ts x_0,n\ts -\ts y_0,t_0\ts -\ts m\ts -\ts n)\in S.
\end{cases}
\end{equation}
As before, we give $\tilde W $ a cell structure, by taking the triangulation inherited from $S$, deleting the horizontal bent edges, 
and deleting the vertices of the horizontal bent edges that are in $\partial \tilde W$.
This might create some non-straight edges in $\partial \tilde W$ as illustrated in Example \ref{extet}.

{
\psfrag{c}{$\Cmn$}
\psfrag{s}{$S$}
\psfrag{y}{$Y$}
\psfrag{w}{$\tilde W$}
\begin{gather*}
\begin{matrix}\epsfig{file=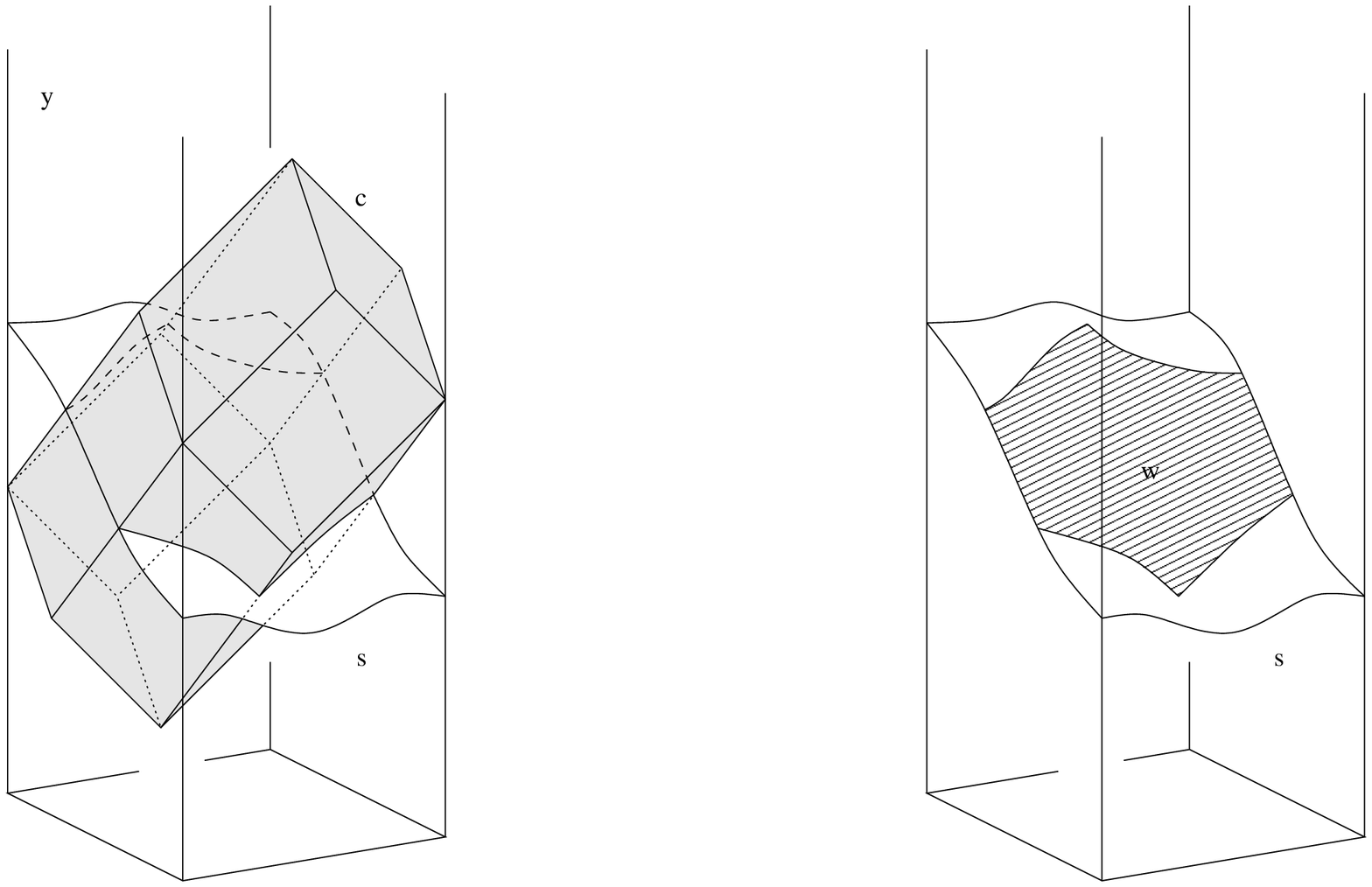,height=5cm}\end{matrix}
\\
\nonumber\text{The complex $\tilde W$.}
\end{gather*}
}
To reduce the number of cases to treat later (in particular in the definition of the constant $\epsilon$), 
we replace these boundary edges \put(3,-2){\epsfig{file=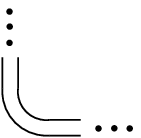,height=.6cm}}\hspace{.8cm} 
by the corresponding straight edges \put(3,-2){\epsfig{file=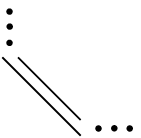,height=.6cm}}\hspace{.8cm}. 
This operation does not affect the combinatorics of the cell structure.
It's this slightly modified complex that we now call $W$.
Note that $W$ is not a subcomplex of $S$ anymore, but that $W \cap \partial Y$ is still contained in $\partial S$.
As before, we shall often represent $W$ via its projection to $\R^2$.

\begin{Example} \label{eg:W}
For $S$ as in Example \ref{extet}, the projection of $W$ on $\R^2$ now looks like this:
\[
W=\,\begin{matrix}\epsfig{file=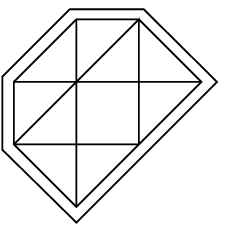,height=1.5cm}\end{matrix}
\]
\end{Example}

Given a matching $M$ of $W$, we define the corresponding matching monomial\dontshow{monomial}
\begin{equation} \label{monomial}
\mu(M) = \prod_{(x,y,t)\in W\cap\tubelattice} f(x,y,t)^{k(x,y,t)} ,
\end{equation}
where $k(x,y,t)$ is given by\dontshow{eme}
\begin{equation}\label{eme}
\begin{split}
k(x,y,t)=1/2\big(\,&\#\,\text{of solid edges incident to}\, (x,y,t) \\
-\:&\#\,\text{of dotted edges incident to}\, (x,y,t)\big) +\epsilon(x,y,t),
\end{split}
\end{equation}
and $\epsilon(x,y,t)\in\{-1,-1/2,0,1/2,1\}$ is a constant depending on the local geometry of $S$ and $\mathcal C$ around $(x,y,t)$ and that we will 
define in\footnote{We believe that the definition of $\epsilon$ is not really needed in order to understand Theorem \ref{tip} or its implications.
The reader who disagrees with that statement is encouraged to read Section \ref{s:tpf} before continuing.
} Section \ref{s:tpf}.
Let $I$ be $W\cap \partial S$ minus its set of isolated points.
The complement of $I$ inside $\partial S$ then admits a disjoint union decomposition
\[
\partial S\setminus I=(\partial S)_+\sqcup (\partial S)_-
\]
where $(\partial S)_+$ and $(\partial S)_-$ are the relative interiors of the closures of the sets
$\big\{(x,y,t)\in\partial S\setminus W\,\big|\,\exists t'< t,\,(x,y,t')\in \Cmn\big\}$
and
$\big\{(x,y,t)\in\partial S\setminus W\,\big|\,\exists t'> t,\,(x,y,t')\in \Cmn\big\}$
respectively.
%\begin{gather*}
%\overline{\big\{(x,y,t)\in\partial S\setminus W\,\big|\,\exists t'< t,\,(x,y,t')\in \Cmn\big\}}\\
%\text{and}\qquad
%\overline{\big\{(x,y,t)\in\partial S\setminus W\,\big|\,\exists t'> t,\,(x,y,t')\in \Cmn\big\}}\qquad\phantom{\text{and}}
%\end{gather*}
%respectively.
{
\psfrag{w}{$W$}
\psfrag{i}{$\scriptstyle I$}
\psfrag{s+}{$\scriptstyle (\partial S)_+$}
\psfrag{s-}{$\scriptstyle (\partial S)_-$}
\begin{gather*}
\begin{matrix}\epsfig{file=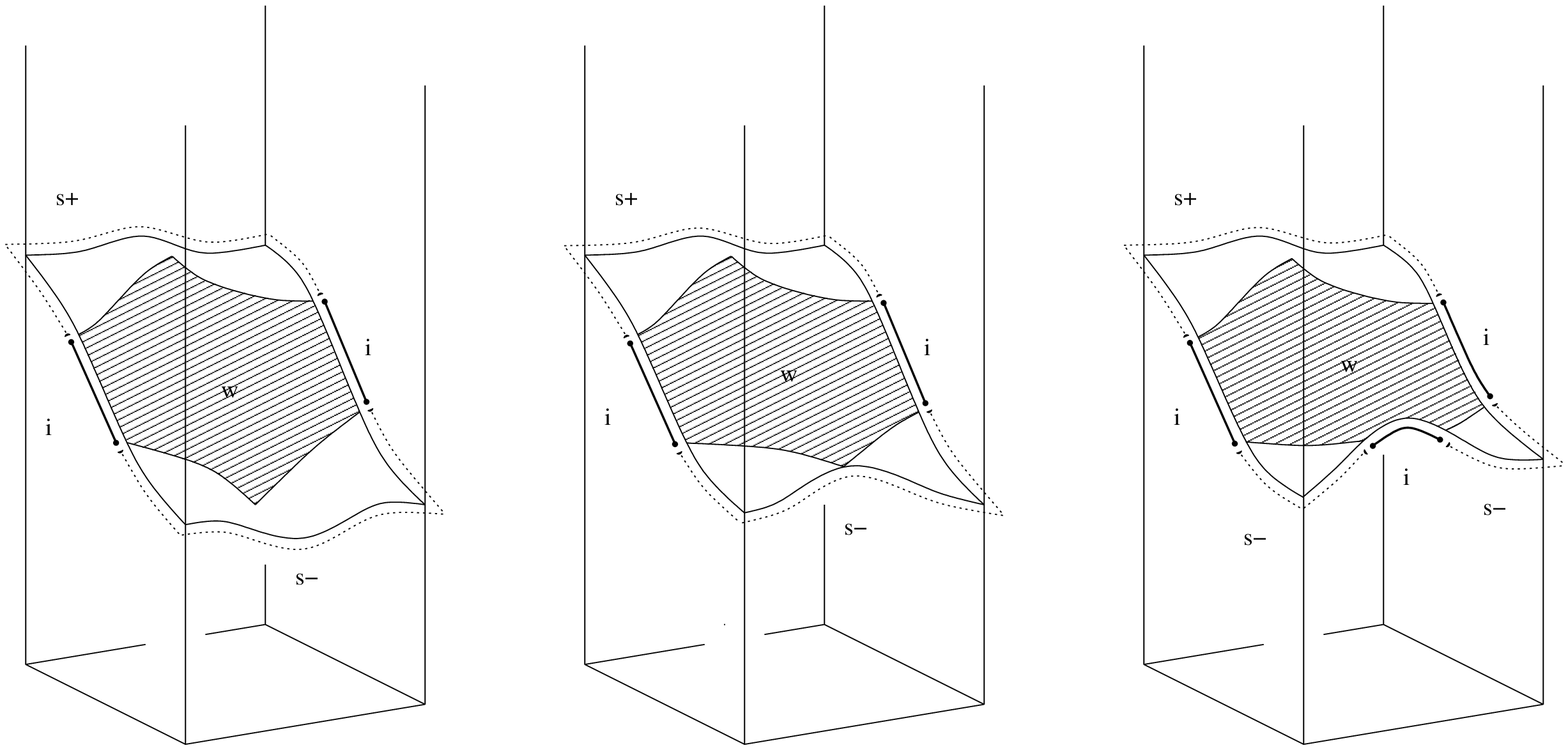,height=5cm}\end{matrix}
\\
\nonumber\text{Some examples of the sets $I$, $(\partial S)_+$, and $(\partial S)_-$.}
\end{gather*}
}

Let $c$ be the Laurent monomial given by\dontshow{imi}
\begin{equation}\label{imi}
c=\prod_{\substack{\text{local maxima}\\ \text{of $(\partial S)_-$}}}\!f(x,y,t)\,\cdot\prod_{\substack{\text{local minima}\\ \text{of $(\partial S)_-$}}}\!f(x,y,t)^{-1},
\end{equation}
where ``local maxima'' and ``local minima'' refers to the time-coordinate.
%$(\partial S)_-$ denotes the part\footnote{A more precise definition of $(\partial S)_-$ is given in (\ref{mpd}).} of $\partial S\setminus \partial W$ that lies below $\mathcal C$.
We then have:\dontshow{tip}
\begin{Theorem}\label{tip}
Let $S$ be a section, $f$ a state, and $(x_0,y_0,t_0)\in \overset{\circ}{Y}\cap\tubelattice$ a point in the future of $S$.
Assume that the light cone $\Cmn$ meets $S$ in at least one point and let $c$, $\mu(M)$ be defined as above. 

Then the value $f(x_0,y_0,t_0)$ given by the bounded octahedron recurrence is determined from $f|_S$ by the formula
\begin{equation*}%\label{eqnmainthm}
f(x_0,y_0,t_0)=c\cdot\sum_{M \in \text{\rm matchings of } W}\mu(M),
\end{equation*}
where $W$, $\mu$, and $c$ are defined in the text above Example \ref{eg:W}, in (\ref{monomial}), and in (\ref{imi}), respectively.
\end{Theorem}

As a corollary, we get the following proof of Theorem \ref{pPt}: 

\begin{proof}[Proof of Theorem \ref{pPt}]
Pick a point $(x_0,y_0,t_0)\in \overset{\circ}{Y}\cap\cL$ and
let $S$ be a section containing $(m-x_0,n-y_0,t_0-m-n)$.
In that case, $W$ consists of a single point and we have $(\partial S)_-=\partial S$. 
There is exactly one matching (the empty matching) and the corresponding matching monomial is $\mu(M)=f(m-x_0,n-y_0,t_0-m-n)$. 
So by Theorem \ref{tip} we get\dontshow{efw}
\begin{equation}\label{efw}
f(x_0,y_0,t_0)=c\cdot f(m-x_0,n-y_0,t_0-m-n)
\end{equation}
as desired.

Now we show that $c$ is independent of $(x_0,y_0,t_0)$ and of $S$.
Given a path $\gamma$ around the boundary of $Y$, by which we mean
a 1-dimensional subcomplex of $\partial Y$ whose projection to $\partial([0,m]\times[0,n])$ is a homeomorphism,
we let $c=c(\gamma)$ be given by
\begin{equation*}
c(\gamma)=\prod_{\substack{\text{local maxima}\\ \text{of $\gamma$}}}\!f(x,y,t)\,\cdot\prod_{\substack{\text{local minima}\\ \text{of $\gamma$}}}\!f(x,y,t)^{-1}.
\end{equation*}
%To see that $\gamma(c)$ is independent of $\gamma$, 
Let
$\gamma'$ be another such path, differing from $\gamma$ by one of the following local moves\dontshow{svd}
\begin{equation}\label{svd}
\begin{matrix}\psfrag{a}{$a$}\psfrag{b}{$b$}\psfrag{e}{$e$}\psfrag{a+b-e}{$ab/e$}\epsfig{file=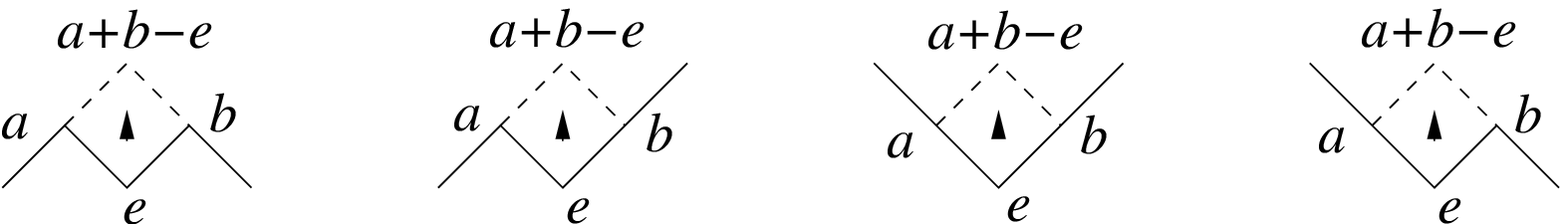,height=1.5cm}\end{matrix}\,\,, 
\end{equation}
where the arrow indicates the direction of increasing time.
We then compute\dontshow{tcf}
\begin{equation}\label{tcf}
\begin{split}
\text{Case 1)}\qquad c(\gamma')&=\cdots(ab/e)\cdots=\cdots ae^{-1}b\cdots=c(\gamma)\\
\text{Case 2)}\qquad c(\gamma')&=\cdots(ab/e)b^{-1}\cdots=\cdots ae^{-1}\cdots=c(\gamma)\\
\text{Case 3)}\qquad c(\gamma')&=\cdots a^{-1}(ab/e)b^{-1}\cdots=\cdots e^{-1}\cdots=c(\gamma)\\
\text{Case 4)}\qquad c(\gamma')&=\cdots a^{-1}(ab/e)\cdots=\cdots e^{-1}b\cdots=c(\gamma).
\end{split}
\end{equation}
Since any two paths can be joined by a sequence of moves like in (\ref{svd}),
we have shown that $c$ is independent of $\gamma$.
Note that in the above computation, we do not distinguish between the corners of $Y$ and the rest of the boundary
since upon identifying $\partial Y$ with the cylinder $\big(\R/(2m+2n)\big)\times \R$, 
all the boundary cases of $(\ref{octrec})$ look the same.

By suitably picking $\gamma$, this computation also proves (\ref{efw}) for $(x_0,y_0,t_0)\in\partial Y$:
let $\gamma$ be the union of the two shortest geodesics in $\partial Y$ between $(x_0,y_0,t_0)$ and $(m-x_0,n-y_0,t_0-m-n)$.
They both have slope 1, so this indeed defines a subcomplex of $\partial Y$.
\[
\psfrag{d}{$\scriptstyle(m-x_0,n-y_0,t_0-m-n)$}
\psfrag{(x,t)}{$\scriptstyle(x_0,y_0,t_0)$}\psfrag{g}{$\gamma$}
\epsfig{file=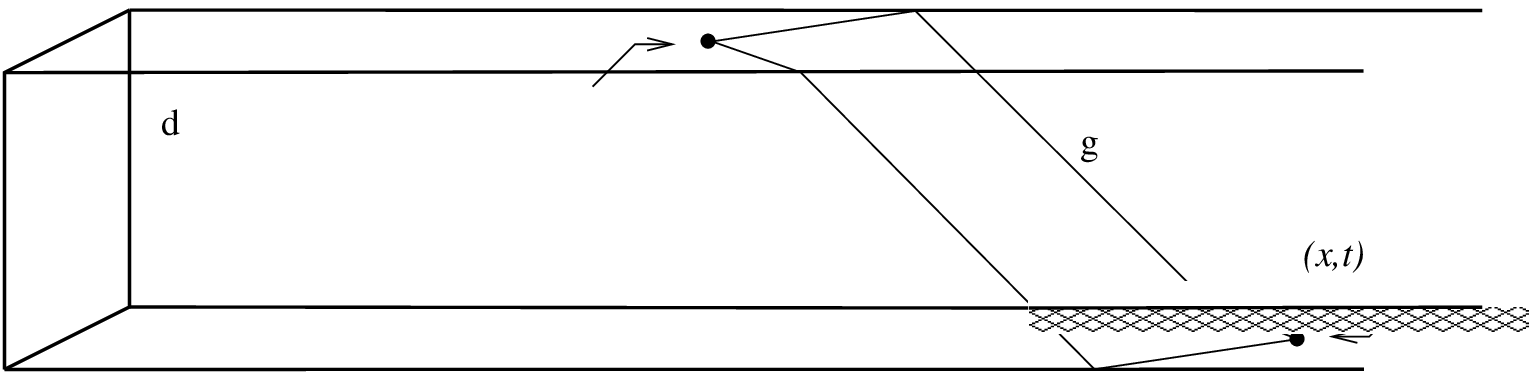,height=1.8cm}
\]
\noindent The path $\gamma$ has exactly two extrema: one maximum $(x_0,y_0,t_0)$ and one minimum $(m-x_0,n-y_0,t_0-m-n)$. So by definition 
$c=f(x_0,y_0,t_0)f(m-x_0,n-y_0,t_0-m-n)^{-1}$, which is exactly (\ref{efw}).
\end{proof}

\begin{Remark}
If $S$ doesn't meet $\Cmn$ then, as it is stated, Theorem \ref{tip} doesn't allow us to compute $ f(x_0,y_0,t_0) $ directly.  
But combining Theorems \ref{tip} and  \ref{pPt},
we may calculate $ f(x_0,y_0,t_0) $ in a non-recursive manner for any point in the future or past of $ S $.
\end{Remark}

\begin{Remark}\label{Trk}
We may\dontshow{Trk} consider variants of the bounded octahedron recurrence where the space coordinates form a strip, or a half-strip, or a quadrant, or a half-plane.
%Letting\dontshow{Trk} $Y=[0,m]\times \R_{\ge 0}\times\R$ or $[0,m]\times \R\times\R$ or $\R_{\ge 0}\times\R\times\R$ or $\R_{\ge 0}^2\times\R$,
%we can consider the obvious generalizations the bounded octahedron recurrence the live on
%$\cL=\{(x,y,t)\in\Z^3\cap Y\,|x+y+t\text{ is even}\}$.
The definitions and results of the above section then extend to these new situations, and can be deduced by a straightforward limit argument.
\end{Remark}

\subsection{The precise formulas}\label{s:tpf}

In this section, we provide the definition of the constant $\epsilon$ used in (\ref{eme}), thus completing the statement of Theorem \ref{tip}.
It is given in the following table:\dontshow{table}

\begin{equation}
\label{table}
\begin{tabular}{|l||c|c|c|c|}\hline
$\epsilon(x,y,t)=$     &       $\hspace{.04cm} t_1>t<t_2\hspace{.04cm}$ & 
                               $\hspace{.04cm} t_1<t<t_2\hspace{.04cm}$ &
                               $\hspace{.04cm} t_1>t>t_2\hspace{.04cm}$ & 
                               $\hspace{.04cm} t_1<t>t_2\hspace{.04cm}$ \\ \hline\hline
\phantom{\Big|}$(x,y,t)\in       \overset{\circ}\Cmn               $ & \multicolumn{4}{c|} {$-1\phantom{-}$} \\ \hline
\phantom{\Big|}$(x,y,t)\in       \partial^{(1)}\Cmn                $ & \multicolumn{4}{c|} {  0            } \\ \hline
\phantom{\Big|}$(x,y,t)\in       \partial^{(2)}\Cmn\setminus I     $ & \multicolumn{4}{c|} { 1/2           } \\ \hline
\phantom{\Big|}$(x,y,t)\in       \partial^{(3)}\Cmn                $ & \multicolumn{4}{c|} {  1            } \\ \hline
\phantom{\Big|}$(x,y,t)\in I^{\circ}\setminus\partial^{\times}\Cmn $ & $-1/2\phantom{-}$ &  0  &  -  & 1/2   \\ \hline
\phantom{\Big|}$(x,y,t)\in \partial^{+}I                           $ &    0              & 1/2 &  0  & 1/2   \\ \hline
\phantom{\Big|}$(x,y,t)\in \partial^{-}I                           $ & $-1/2\phantom{-}$ &  0  & 1/2 &  1    \\ \hline
\phantom{\Big|}$(x,y,t)\in       \partial^{\times}\Cmn             $ &    -              & 1/2 & -   &  -    \\ \hline
\end{tabular}
\end{equation}

\noindent Here, $\overset{\circ}\Cmn$ is the interior of $\Cmn$. We use the notations $\partial^{(1)}\Cmn$, 
$\partial^{(2)}\Cmn$, $\partial^{(3)}\Cmn$, 
$\partial^{\circ}\Cmn$, $\partial^{+}\Cmn$, $\partial^{-}\Cmn$,
and $\partial^{\times}\Cmn$, for the various subsets of its boundary depicted in the following figure:

\begin{center}
\begin{minipage}{6.4cm}
   \epsfig{file=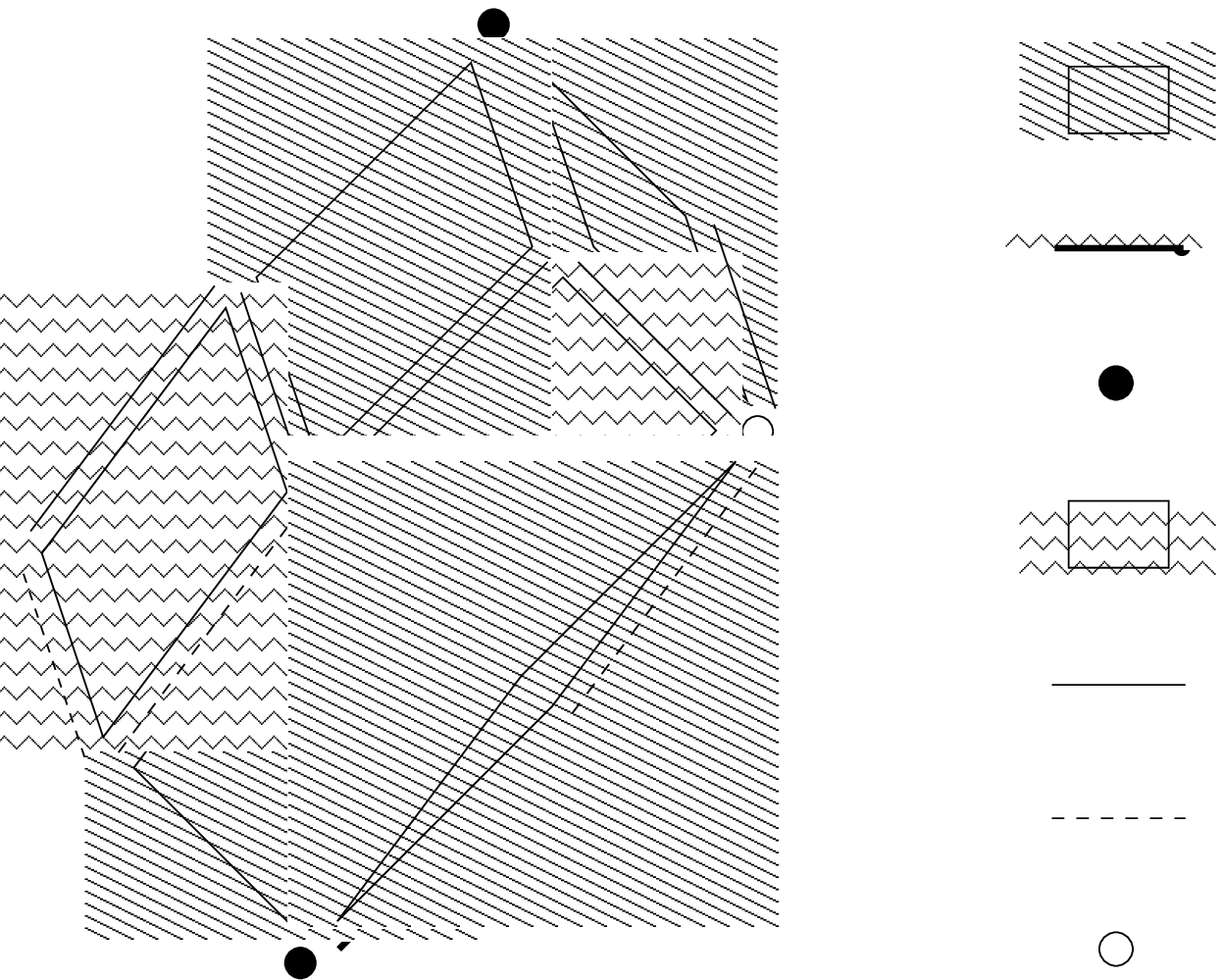,height=5cm}
\end{minipage}
\begin{minipage}{1cm}
\vspace{.3cm}

$\partial^{(1)}\Cmn$
\vspace{.35cm}

$\partial^{(2)}\Cmn$
\vspace{.25cm}

$\partial^{(3)}\Cmn$
\vspace{.35cm}

$\partial^{\circ}\Cmn$
\vspace{.35cm}

$\partial^{+} \Cmn$
\vspace{.25cm}

$\partial^{-} \Cmn$
\vspace{.3cm}

$\partial^{\times}\Cmn$
\end{minipage}
\end{center}
\vspace{.3cm}

As explained before,
the set $I$ is $W\cap\partial Y$ minus its set of isolated points (which necessarily belong to $\partial^{(2)}\Cmn$).
%equivalently, it's the closure of $W\cap (\partial^{\circ}\Cmn \cup \partial^{+} \Cmn \cup \partial^{-} \Cmn)$. 
Its relative interior is denoted 
$I^{\circ}$. Its relative boundary $\partial I$ decomposes into two parts $\partial^+I:=I\cap\overline{(\partial S)_+}$ and $\partial^-I:=I\cap\overline{(\partial S)_-}$.

%, where 
%$\partial^+I\subset \overline{\partial^{+} \Cmn}$ and $\partial^-I\subset \overline{\partial^{-} \Cmn}$.
%Note that $\partial I$ never meets $\partial^\times \Cmn$, and so the above description of $\partial^+I$ and $\partial^-I$ is unambiguous.

If $(x,y,t)\in \partial S$, the numbers $t_1$ and $t_2$ are the heights of its two neighbors $(x_1,y_1,t_1)$ and $(x_2,y_2,t_2)$
in $\partial S$. If both $(x_1,y_1,t_1)$ and $(x_2,y_2,t_2)$ belong to $W$, 
we order them so that $t_1\le t_2$. Otherwise, we let $(x_1,y_1,t_1)\in W$ and $(x_2,y_2,t_2)\not \in W$.

\section{Proof of theorem \ref{tip}}\label{the big proof}

The\dontshow{the big proof} proof follows the first of the two proofs of Theorem \ref{thmoct} presented in \cite{Spe04}. 
It will be in three steps. 
First, we reduce ourselves to the case when $S$ is contained in $\Cmn$.
Secondly, if $S\subset \Cmn$,
we build an auxiliary complex $\bW$ whose matchings are in bijection with those of $W$, but where the formula for the 
matching monomials is simpler. Finally we use induction on the distance between $S$ and $(x_0,y_0,t_0)$ to prove our auxiliary formula.

Note that if $m=1$ or $n=1$ then $\overset{\circ}{Y}\cap\tubelattice$ is empty, and hence Theorem \ref{tip} is trivially true.
So from now on let's assume that $m,n\ge 2$.

\subsection{First Step}

\begin{Lemma}\label{lmun}
Assume\dontshow{lmun} that Theorem \ref{tip} holds for all $S\subset \Cmn$. Then it holds for all $S$.
\end{Lemma}

\begin{proof}
We prove this by induction on the volume between $S$ and $\Cmn$. Let $S$ be a section and let $S'$ be another section, 
just a little bit closer to $\Cmn$, and agreeing with $S$ inside $\Cmn$. 
More precisely, we assume that the volume between $S$ and $S'$ consists of a single 3-cell, not contained in $\Cmn$.

Let $W$, $W'$ be the complexes constructed in section \ref{sec:Tbc} corresponding to $S$, $S'$ respectively, and let $c$, $c'$ be the associated constants (\ref{imi}).
By induction, we assume that 
$f(x_0,y_0,t_0)= c'\cdot\sum\mu(M')$, and we want to show that $f(x_0,y_0,t_0)= c\cdot\sum\mu(M)$. 
Here $M$ runs over all matchings of $W$, and $M'$ runs over all matchings of $W'$.
It will therefore be enough to show that 
$c'\cdot\sum\mu(M')=c\cdot\sum\mu(M)$.

Since the difference between $S$ and $S'$ occurs only outside of $\Cmn$, we have $W'=W$ 
and to every matching $M$ of $W$ there is a corresponding matching $M'$ of $W'$.
However, we don't always have $\mu(M')=\mu(M)$ since
the definition (\ref{monomial}) of $\mu(M)$ involves the coefficients $\epsilon$ given in (\ref{table}). Those depend on the
height $t_2$ which might change between $S$ and $S'$, thus affecting the sum.
We will prove our claim by showing that in each case $c\cdot\mu(M)=c'\cdot\mu(M')$.

Case 1) The move between $S$ and $S'$ happens above $\Cmn$. We check in row 6 of (\ref{table}) that $\epsilon$ doesn't
depend on $t_2$. Therefore $\mu(M')=\mu(M)$ for each matching $M$. The constant $c$ only depends on things happening below $\Cmn$, thus
$c'=c$, which proves our claim. \qe

Case 2) The move happens below $\Cmn$. We can again distinguish between various cases.

\noindent Case 2.1) This happens far enough from $\Cmn$, so that $\mu(M')=\mu(M)$. Then we only need to show that $c'=c$.
The only moves that affect the value of $c$ are the ones near the boundary $\partial Y$. 
Restricted to $(\partial S)_-$, these moves look exactly like (\ref{svd}).
The verification that $c'=c$ is identical to (\ref{tcf}). \qe

If the move happens close to $\Cmn$, then either of the points labeled $a$ or $b$ in (\ref{svd})
could fail to be in $(\partial S)_-$. They would
therefore have to be in $W$.

\noindent Case 2.2) Suppose $a\not\in (\partial S)_-$ and $b\in (\partial S)_-$. Then we have
\begin{equation*}\begin{split}
\text{Case 2.2.1)}\qquad c'&=(ab/e)\cdots=a\cdot(e^{-1}b\cdots)=ac\\
\text{Case 2.2.2)}\qquad c'&=(ab/e)b^{-1}\cdots=a\cdot(e^{-1}\cdots)=ac\\
\text{Case 2.2.3)}\qquad c'&=(ab/e)b^{-1}\cdots=a\cdot(e^{-1}\cdots)=ac\\
\text{Case 2.2.4)}\qquad c'&=(ab/e)\cdots=a\cdot(e^{-1}b\cdots)=ac.\\
\end{split}\end{equation*}

\noindent 
In all cases we have $c'=ac$, so we need to show that 
$\mu(M')=a^{-1}\cdot\mu(M)$. 
Indeed the
only difference between $\mu(M')$ and $\mu(M')$ is the constant $\epsilon$ in row 7 of table (\ref{table}) that influences the exponent of $a$.
Let us call $\epsilon'$ and
$\epsilon$ the constants used in $\mu(M')$ and $\mu(M)$ respectively. For $\epsilon'$ we have $t'_2>t$, and for $\epsilon$ we have
$t_2<t$ (here, we use the notation of table (\ref{table}), namely $t$ is the height of the 
point $a$ and $t_2$ the height of the point $e$). Regardless of the value of $t_1$, we have $\epsilon'=\epsilon-1$ which shows that
$\mu(M')=a^{-1}\cdot\mu(M)$. \qe

\noindent Case 2.3) If $b\not\in (\partial S)_-$ and $a\in (\partial S)_-$ then the situation is symmetric to that in 2.2). \qe

\noindent Case 2.4) If neither $a$ nor $b$ are in $(\partial S)_-$ then an argument similar to 2.2) shows that $c'=abc$ and that 
$\mu(M')=a^{-1}b^{-1}\cdot\mu(M)$. \qe
%
%We have shown that in all cases $c'\cdot\mu(M')=c\cdot\mu(M)$, which finishes our inductive proof.
\end{proof}

\subsection{Second Step}\label{s:2St}

From now on, we assume that $S\subset \Cmn$. Following\dontshow{s:2St} \cite{Spe04}, we introduce a variant $\bW$ of $W$.
Instead of starting with $S\cap\overset{\circ}{\Cmn}$, we take $S$ itself and remove the horizontal bent edges (see section \ref{3p1}): that's $\bW$.
The part of $\bW$ that is outside of $W$ corresponds to the part of $S$ that lies in the upper or lower boundary of $\Cmn$.
It comes with the following simple triangulation:

\begin{equation}\label{pattern}\begin{matrix}
\psfrag{x}{$x_0$}\psfrag{y}{$y_0$}
\psfrag{0}{$0$}\psfrag{m}{$m$}\psfrag{n}{$n$}
\epsfig{file=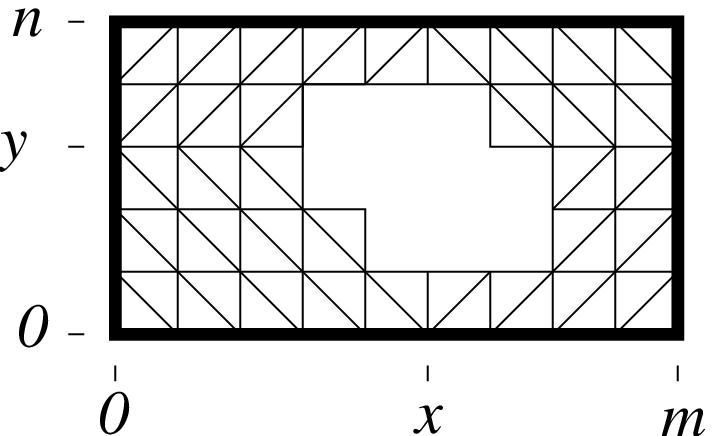,height=2.5cm}\end{matrix}
\end{equation}

The edges of $\partial W$ that are diagonals in the 2-dimensional projection, correspond to horizontal bent edges in $S$.
As such, they get removed from the cell structure of $\bW$.
This gives us a way of describing $\bW$ from $W$ entirely in terms of the 2-dimensional projection.
Namely, $\bW$ is obtained from $W$ by removing all diagonal edges of $\partial W$, and then completing with the pattern (\ref{pattern}).

\begin{Example}
Let $S$ be the section obtained from the section of Example \ref{extet} by projecting it down to $\Cmn$. 
Then the corresponding complex $W$ is as shown in Example \ref{eg:W}. 
The auxiliary complex $\bW$ is obtained from $W$ by completing with the  pattern (\ref{pattern}):

\begin{equation*}
\psfrag{s1}{$W=$}\psfrag{s2}{$\bW=$}
\begin{matrix}\epsfig{file=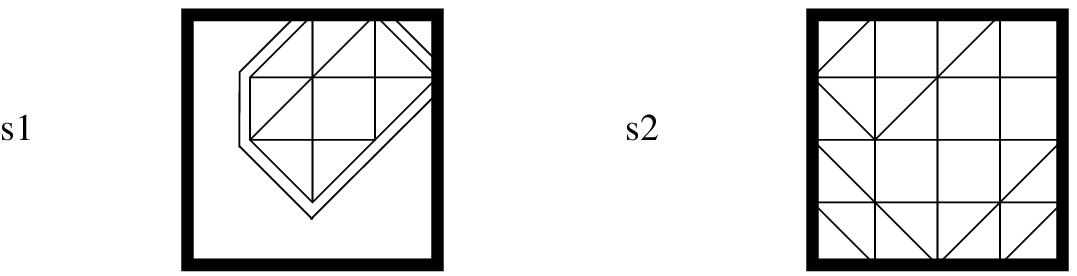,height=1.8cm}\end{matrix}
\end{equation*}
\end{Example}

\noindent 
As before, given a matching $M$ of $\bW$, we define a corresponding matching monomial\dontshow{mmp}

\begin{equation}\label{mmp}
\overline\mu(M)=\prod_{(x,y,t)\in \bW\cap\tubelattice}f(x,y,t)^{\overline k(x,y,t)},
\end{equation}
where
\[
\begin{split}
\overline k(x,y,t)=1/&2\big(\,\#\,\text{of solid edges incident to}\, (x,y,t) \\
&-\:\#\,\text{of dotted edges incident to}\, (x,y,t)\big) +\overline \epsilon(x,y,t),
\end{split}
\]
and the value of $\overline \epsilon(x,y,t)$ is given by:\dontshow{tbp}

\begin{equation}\label{tbp}
\begin{tabular}{|l||c|c|c|}\hline
$\overline\epsilon(x,y,t)=$     & $\hspace{.2cm}  t_1>t<t_2\hspace{.2cm} $ & 
                                  $\hspace{.2cm}  t_1<t<t_2\hspace{.2cm} $ &
                                  $\hspace{.2cm}  t_1<t>t_2\hspace{.2cm} $ \\ \hline\hline
\phantom{\Big|}$(x,y,t)\in S\cap \overset{\circ}Y$          & \multicolumn{3}{c|} {$-1\phantom{-}$} \\ \hline
\phantom{\Big|}$(x,y,t)\in S\cap \partial^{(1)}Y $          &  $-1/2$  &  0  & 1/2   \\ \hline
\phantom{\Big|}$(x,y,t)\in S\cap \partial^{(2)}Y $          &  0  &  1/2  & 1  \\ \hline
\end{tabular}
\end{equation}

\noindent Here $\partial^{(2)}Y$ denotes the four vertical edges of $Y$, and $\partial^{(1)}Y=\partial Y\setminus \partial^{(2)}Y$ 
denotes the interior of the facets of $Y$. As in (\ref{table}), the numbers $t_1$ and $t_2$ refer to the heights of the neighbors
of $(x,y,t)$ in $\partial S$. We order them so that $t_1\le t_2$.

\begin{Lemma} \label{lem9}
There\dontshow{lem9} is a natural bijection between matchings of $W$ and matchings of $\bW$. 
Under that bijection we have $c\cdot\mu(M)=\overline\mu(M)$, where $c$, $\mu(M)$, and $\overline\mu(M)$ are defined in 
(\ref{imi}), (\ref{monomial}) and (\ref{mmp}) respectively.
\end{Lemma}

\begin{proof}
The pattern (\ref{pattern}) used to complete $W$ is composed of four ``Young diagrams'' (\ref{TYD}.i), one in each corner.\dontshow{TYD}
%\break\put(0,-5){\epsfig{file=completion.eps,height=1cm}}\hspace{1.30cm}
It is easy to see, starting from the triangle in the corner, that the only possible matching of this region is (\ref{TYD}.ii). 
\begin{equation}\label{TYD}
\text{(i)}\begin{matrix}
\epsfig{file=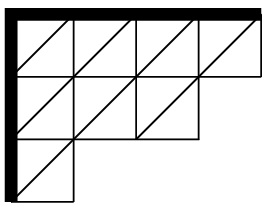,height=1cm}
\end{matrix}\qquad\qquad
\text{(ii)}\begin{matrix}
\epsfig{file=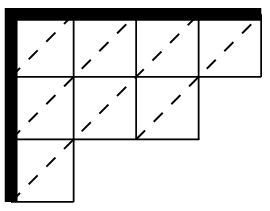,height=1cm}
\end{matrix}
\end{equation}
%\put(4,-5){\epsfig{file=completion2.eps,height=1cm}}\hspace{1.45cm}. 
The remaining edges of $\bW$ are exactly the internal edges of $W$. 
The matchings of $\bW$ and of $W$ are thus naturally in bijection.

Let us rewrite (\ref{imi}) as
\begin{equation*}
c=\prod_{(x,y,t)\in (\partial S)_-}f(x,y,t)^{j(x,y,t)},
\end{equation*}
where\dontshow{defj}
\begin{equation}\label{defj}
j=j(x,y,t)=\begin{cases}
\;1&\text{ if $(x,y,t)$ is a local maximum in $(\partial S)_-$,}\\ 
-1&\text{ if $(x,y,t)$ is a local minimum in $(\partial S)_-$,}\\ 
\;0&\text{ if $(x,y,t)$ is not a local extremum in $(\partial S)_-$.}\\ 
\end{cases}
\end{equation}
In order to show that $c\cdot\mu(M)=\overline \mu(M)$, we need to check at each vertex $(x,y,t)$ that
\begin{equation}\label{aux}
\left[\parbox{5.15cm}{
$1/2\cdot\big(\,\#$ of solid edges in $W$

$\phantom{*}-\:\#$ of dotted edges in $W\big)$ $+\,\epsilon$}\right]+j
=
\left[\parbox{5.2cm}{
$1/2\cdot\big(\,\#$ of solid edges in $\bW$

$\phantom{*}-\:\#$ of dotted edges in $\bW\big)$ $+\,\overline\epsilon$}\right],
\end{equation}
where $\epsilon$, $j$ and $\overline \epsilon$ are defined in (\ref{table}), (\ref{defj}), and (\ref{tbp}) respectively.
For points that are not in $(\partial S)_-$, we set $j=0$.

This is done by a case by case study.

Case 1) The point $(x,y,t)$ is in the interior of $W$. Then the number of solid and dotted edges
incident to our point doesn't change. 
We have $\epsilon=\overline\epsilon=-1$ and $j=0$, therefore (\ref{aux}) holds. \qe

Case 2) The point $(x,y,t)$ is not in $W$. Then we need to show that 
$j=1/2\cdot(\#$ of solid edges in $\bW$
$-\:\#$ of dotted edges in $\bW)$ $+\,\overline\epsilon$.

\noindent Case 2.1) The point $(x,y,t)$ is in the interior of $Y$. Then $j=0$ and $\overline \epsilon=-1$. 
The local picture looks like this
\put(5,-2){\epsfig{file=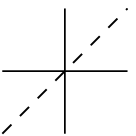,height=.7cm}}\hspace{.9cm} or like this
\put(5,-2){\epsfig{file=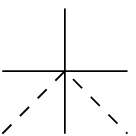,height=.7cm}}\hspace{.9cm}. In both cases, there are 4 solid and 2 dotted edges incident to our point and
we can check that (\ref{aux}) holds. \qe

\noindent Case 2.2) The point $(x,y,t)$ lies in $\partial^{(1)}Y$. 

\noindent Case 2.2.1) If $x\not =x_0$ and $y\not =y_0$, then our point is not a local extremum in $\partial S$ and therefore
$j=\overline \epsilon = 0$. The local picture is 
\put(2,-1){\epsfig{file=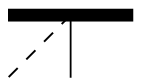,height=.45cm}}\hspace{1cm}, where the thick line represents the boundary $\partial S$.
There's 1 solid and 1 dotted edge: (\ref{aux}) holds. \qe

\noindent Case 2.2.2) If either $x=x_0$ or $y=y_0$, then $(x,y,t)\in\partial^{(2)}\Cmn\cap\partial Y$ and it is a local extremum
in $\partial S$. The local picture is \put(2,-1){\epsfig{file=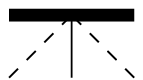,height=.45cm}}\hspace{1cm}.
If it's a maximum, then $(x,y,t)$ lies in the top part of $\partial^{(2)}\Cmn$, and therefore $(x,y,t)\not\in (\partial S)_-$.
We have $j=0$ and $\overline \epsilon=1/2$, and one checks (\ref{aux}).
If it's a minimum, then $(x,y,t)\in (\partial S)_-$, $j=-1$, $\overline \epsilon=-1/2$, and one checks (\ref{aux}). \qe

\noindent Case 2.3) The point $(x,y,t)$ lies in $\partial^{(2)}Y$. 
The local picture is \put(4,-1){\epsfig{file=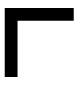,height=.45cm}}\hspace{.65cm}. The point $(x,y,t)$ is either
a maximum or a minimum. If it's a minimum, then $(x,y,t)\not\in (\partial S)_-$, $j=\overline \epsilon=0$, and (\ref{aux}) holds.
If it's a maximum, then then $(x,y,t)\in (\partial S)_-$, $j=\overline\epsilon=1$, and again (\ref{aux}) holds. \qe

Case 3) The point $(x,y,t)$ is on the boundary $\partial W$ but not in $\partial Y$. 
Then $j=0$. We distinguish between the following cases:

\noindent Case 3.1) If it lies in $\partial^{(1)}\Cmn$, then $\epsilon=0$ and $\overline \epsilon=-1$.
The possible local situations in $W$ are listed below,
next to the corresponding $\bW$ pictures. The double line represents $\partial W$.
\[
\epsfig{file=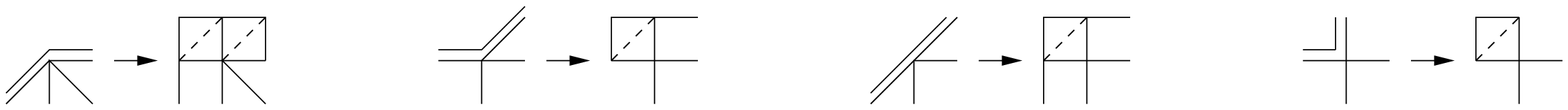,width=12.5cm}
\]
\noindent In the leftmost case, we added 2 solid and 1 dotted edge, in the three other cases, we just added 2 solid edges.
But in all cases the number $1/2\cdot(\#$ of solid edges $-\:\#$ of dotted edges$)$ got increased by 1 in the process of
replacing $W$ by $\bW$. This is exactly what is needed to make sure (\ref{aux}) holds. \qe

\noindent Case 3.2) If it lies in $\partial^{(2)}\Cmn$, then $\epsilon=1/2$ and $\overline \epsilon=-1$.
We list the possible local situations in $W$ and in $\bW$:
\[
\epsfig{file=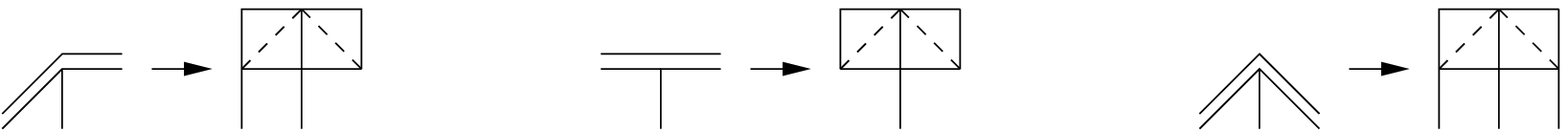,width=8.8cm}
\]
\noindent In each case we added three new solid edges. Therefore $1/2\cdot(\#$ of solid edges $-\:\#$ of dotted edges$)$
got increased by 3/2. This shows (\ref{aux}).

\noindent Case 3.3) If it lies in $\partial^{(3)}\Cmn$, then $\epsilon = 1$ and $\overline \epsilon=-1$. The picture looks like this
\put(3,-3.5){\epsfig{file=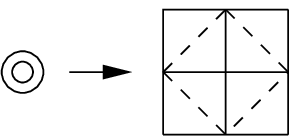,height=.75cm}}\hspace{1.9cm}.
The exponent in $\mu(M)$ is $1/2\cdot(\#$ of solid edges $-\:\#$ of dotted edges$)+\epsilon
=1/2\cdot(0-0)+1=1$.
The exponent in $\mu(\overline M)$ is $1/2\cdot(\#$ of solid edges $-\:\#$ of dotted edges$)+\overline \epsilon
=1/2\cdot(4-0)-1=1$.
Therefore (\ref{aux}) holds. \qe

Case 4) The point $(x,y,t)$ lies in $I$. Then it's not in $(\partial S)_-$ and therefore $j=0$.

\noindent Case 4.1) It lies in $I^\circ$. We can distinguish between 
the following two cases: either it lies in $I^\circ\setminus \partial^\times \Cmn$ or it lies in $\partial^\times \Cmn$,
but in either case, the number of solid and dotted edges
incident to our point doesn't change. So we only need to check that $\epsilon=\overline \epsilon$, which is done by direct 
inspection of (\ref{table}) and (\ref{tbp}). \qe

\noindent Case 4.2) If it lies in $\partial^+I$, then the local situations look like this:
\[
\epsfig{file=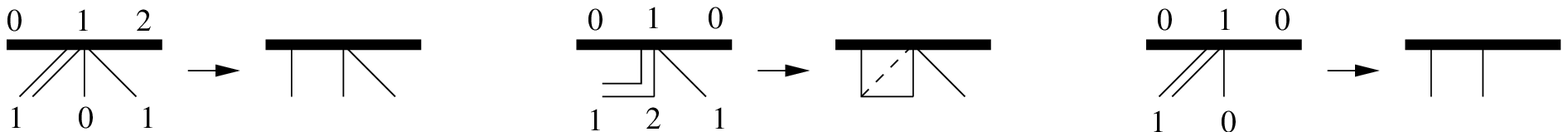,height=1cm}
\]
\noindent Here, we have indicated the heights (up to a constant) of the various points in $S$. As before, the thick line denotes the 
boundary $\partial S$. In each case, the quantity 
$1/2\cdot(\#$ of solid edges $-\:\#$ of dotted edges$)$ 
%$[\frac{1}{2}$ $\#$ internal edges $-\:\#$ matches$]$ 
remained unchanged. So we only
need to check that $\epsilon=\overline \epsilon$. Indeed, in the first case, our point is not an extremum therefore
$\epsilon=\overline \epsilon=0$ and in the other two cases our point is a maximum therefore $\epsilon=\overline \epsilon=1/2$. \qe

\noindent Case 4.3) If it lies in $\partial^-I$, then the local situations look exactly like in 4.2), except that we need to change all
the heights to minus their value. Again the quantity 
$1/2\cdot(\#$ of solid edges $-\:\#$ of dotted edges$)$ 
remains unchanged, and we 
need to check that $\epsilon=\overline \epsilon$. Indeed, one verifies that $\epsilon=\overline \epsilon=0$ in 
the first case and $\epsilon=\overline \epsilon=-1/2$ in the other two cases. \qe

Case 5) The point $(x,y,t)$ lies in $W\cap\partial Y$ but not in $I$. This can only happen if it 
belongs to $\partial^{(2)}\Cmn$, therefore $\epsilon=1/2$, and the local situation
is like this: \break\put(0,-1){\epsfig{file=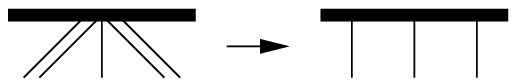,height=.43cm}}\hspace{2.9cm}.
Again, the quantity 
$1/2\cdot(\#$ of solid edges $-\:\#$ of dotted edges$)$ 
remains unchanged, 
so in order to verify (\ref{aux}), we need to check that $\epsilon+j=\overline \epsilon$.

\noindent Case 5.1)
If the point lies above $\partial^\circ\Cmn$, then it's not in $(\partial S)_-$, and therefore $j=0$. Also, it's 
a local maximum in $\partial S$ so 
$\overline \epsilon=1/2$. We then check that $\epsilon+j=\overline \epsilon$. \qe

\noindent Case 5.2)
If the point lies below $\partial^\circ\Cmn$, then it's in $(\partial S)_-$. It's a local minimum in $\partial S$, therefore
$j=-1$ and $\overline \epsilon=-1/2$. We again check that $\epsilon+j=\overline \epsilon$. \qe

This last case completes the proof of Lemma \ref{lem9}.
\end{proof}

\subsection{Third Step}

Recall that $S\subset \mathcal C$ and that want to show\dontshow{qnm}
\begin{equation}\label{qnm}
f(x_0, y_0,z_0) = c\cdot \sum_{M \in \text{matchings of } W} \mu(M).
\end{equation}
If $S$ contains the point $(x_0,y_0,t_0)$, then we have $c=1$, $W=\{(x_0,y_0,t_0)\}$ and $\epsilon(x_0,y_0,t_0)=1$.
There is only one matching and equation (\ref{qnm}) holds.
So in order to prove the theorem, we just need to show that the RHS of (\ref{qnm}) is independent of $S$.
By Lemma \ref{lem9}, we have 
\begin{equation}\label{eqnmainthm'}
 c\cdot \sum_{M \in \text{matchings of } W} \mu(M) = \sum_{{M} \in \text{matchings of } \bW} \overline\mu({M}),
\end{equation}
so it's enough to show that the RHS of (\ref{eqnmainthm'}) is independent of $S$.

Let $ Z(\bW) $ denote the right hand side of (\ref{eqnmainthm'}).  
More generally, given a planar complex $V$ and numbers $\epsilon(v)$ attached to every vertex $v\in V$, 
then to any state $f$ on $V$ we may assign the quantity\dontshow{redefZ}
\begin{equation}\label{redefZ}
Z(V,\epsilon,f):=\sum_{M \in \text{matchings of } V}\quad\prod_{v\in V}f(v)^{k(v,\epsilon)}
\end{equation}
where 
$
k(v,\epsilon)=1/2\cdot(\,\#$ of solid edges incident to $v 
-\:\#$ of dotted edges incident to $v) +\epsilon(v)$.
%\[
%\begin{split}
%\overline k(v)=1/&2\big(\,\#\,\text{of solid edges incident to}\,\, v \\
%&-\:\#\,\text{of dotted edges incident to}\,\, v\big) +\overline \epsilon(v),         [\cite[Proposition 7]{Spe04}]
%\end{split}
%\]
The following lemma is borrowed from Speyer:
\begin{Lemma}[{\cite[Proposition 7]{Spe04}}]\label{lemUR}
Let\dontshow{lemUR} $V$ be a planar complex, and let
$V'$ be the complex obtained from $V$ by adding a double edge between two vertices of a given 2-face.
Then $Z(V,\epsilon,f)=Z(V',\epsilon,f)$.
\end{Lemma}

\begin{proof}
Every matching of $V$ extends in a unique way to a matching of $V'$.
Of the two additional edges, one will be solid and one will be dotted.
The quantity $k(v,\epsilon)$ remains unchanged, and therefore so does $Z(V)$.
\end{proof}

Any two sections can be joined by a sequence $S_1, S_2,\ldots $ such that the volume between $S_i$ and $S_{i+1}$ consists of a single 3-cell,
so we're reduced to proving the following lemma:

\begin{Lemma}\label{LsL}
Let\dontshow{LsL} $S,S'\subset\mathcal C$ be two sections such that the volume between $S$ and $S'$ consists of a single 3-cell,
and let $\bW$, $\bW'$ be the corresponding complexes, as defined in section \ref{s:2St}.
Then
$$
Z(\bW)=Z(\bW').
$$
\end{Lemma}

\begin{proof}
Since $\Cmn$ doesn't contain 1/4-octahedra, the volume between $S$ and $S'$ is either a tetrahedron, an octahedron, or a 1/2 octahedron.

Case 1) If the volume between $S$ and $S'$ is a tetrahedron, 
then $S$ and $S'$ share the same set of vertices.
The projection to $\R^2$ of $S'$ is obtained from that of $S$ by removing the diagonal of a square and replacing it with the other diagonal.
Both these diagonals correspond to horizontal bent edges (see section \ref{3p1}) and are thus absent from $\bW$ and $\bW'$.
It follows that $\bW=\bW'$ and in particular that $Z(\bW)=Z(\bW')$. \qe

Case 2) If the volume between $S$ and $S'$ is an octahedron, then we have the following possible local situations\dontshow{elg}
for $\bW$ and $\bW'$:
\begin{equation} \label{elg}
\begin{matrix}\epsfig{file=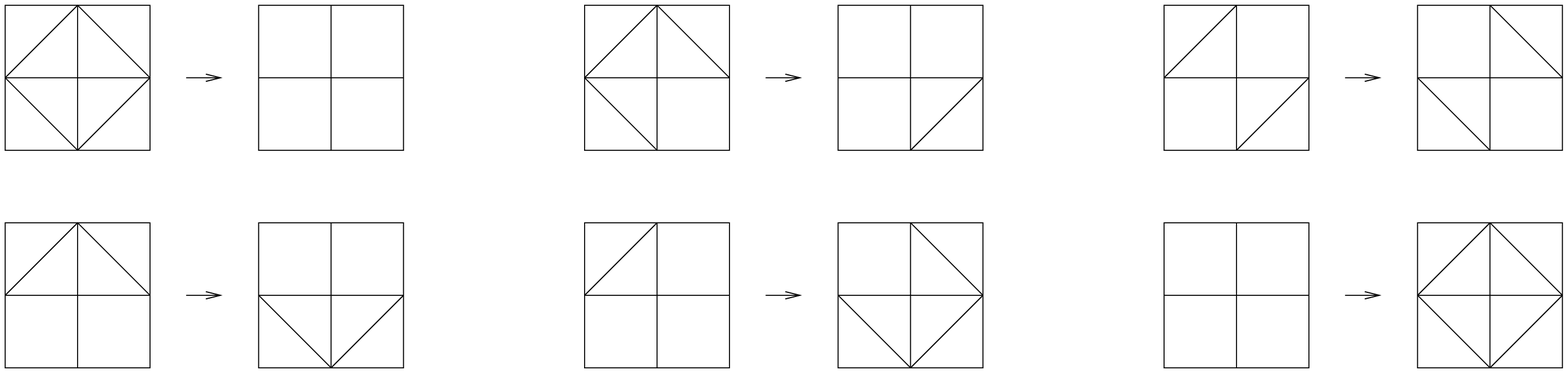,width=9cm}\end{matrix}
\end{equation}
The pictures of $\bW$ are on the left and those of $\bW'$ are on the right. 
We want to show that in each case $Z(S,\overline\epsilon,f)=Z(S',\overline\epsilon',f')$, where $f$ and $f'$ are given respectively by
$$
\begin{matrix}h&a&k\\b&\hspace{.5cm}e\hspace{.5cm}&c\\p&d&q\end{matrix} \hspace{1.5cm}\text{   and   }\hspace{1.5cm}
\begin{matrix}h&a&k\\b&\;\;(ad+bc)/e\;\;&c\\p&d&q,\end{matrix}
$$
and $\overline\epsilon=\overline\epsilon'$.
By Lemma \ref{lemUR}, we may replace (\ref{elg}) by the following complexes without altering the value of $Z$:\dontshow{eq.lg}
\begin{equation}\label{eq.lg}
\begin{matrix}\epsfig{file=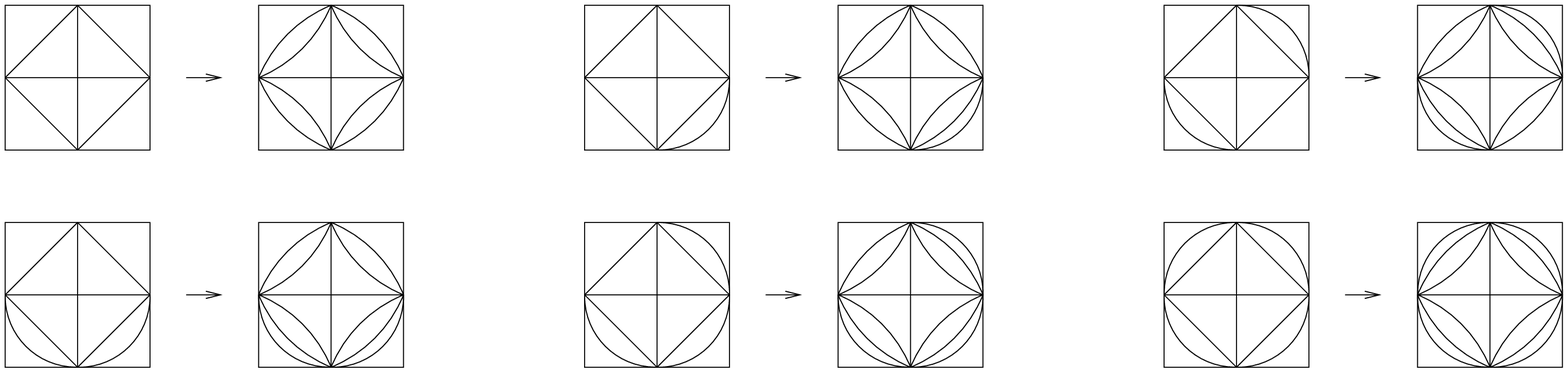,width=9cm}\end{matrix}
\end{equation}
The RHS of (\ref{eq.lg}) is now always obtained from the LHS by the same operation of replacing 
\put(4,-2){\epsfig{file=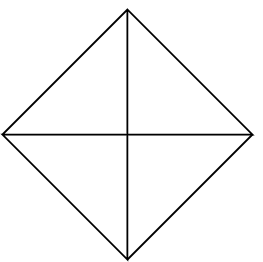,height=.8cm}}\hspace{1cm} by
\put(4,-2){\epsfig{file=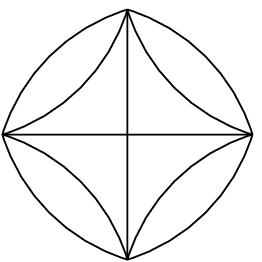,height=.8cm}}\hspace{1cm}.
So we just have to check that\dontshow{Spt}
\begin{equation}\label{Spt}
\psfrag{a}{$a$}
\psfrag{b}{$b$}
\psfrag{c}{$c$}
\psfrag{d}{$d$}
\psfrag{e}{$e$}\psfrag{(ad+bc)/e}{$(ad+bc)/e$}
Z\left(\begin{matrix}\epsfig{file=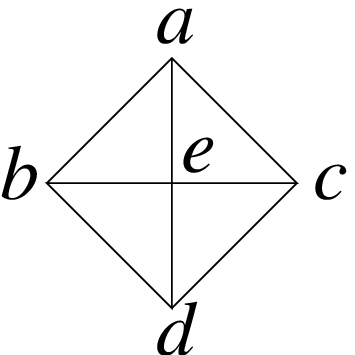,height=1.5cm}\end{matrix}\right)=
Z\left(\begin{matrix}\epsfig{file=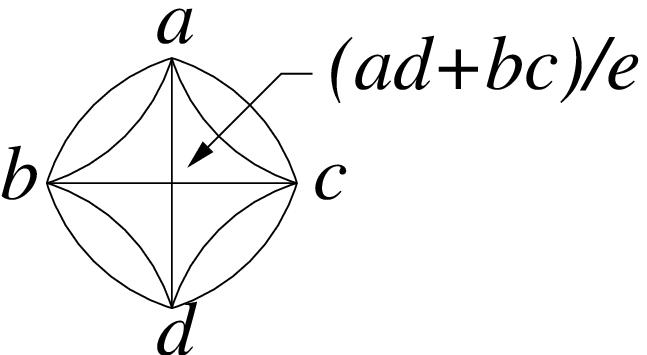,height=1.5cm}\end{matrix}\hspace{.4cm}\right).
\end{equation}
The matchings $M$ of the LHS of (\ref{Spt}) are not exactly in bijection with the matchings $M'$ of the RHS, but there is the 
following correspondence \cite[Section 4.2]{Spe04}:\dontshow{bfk}
\begin{equation} \label{bfk}
\phantom{.}\begin{matrix}\epsfig{file=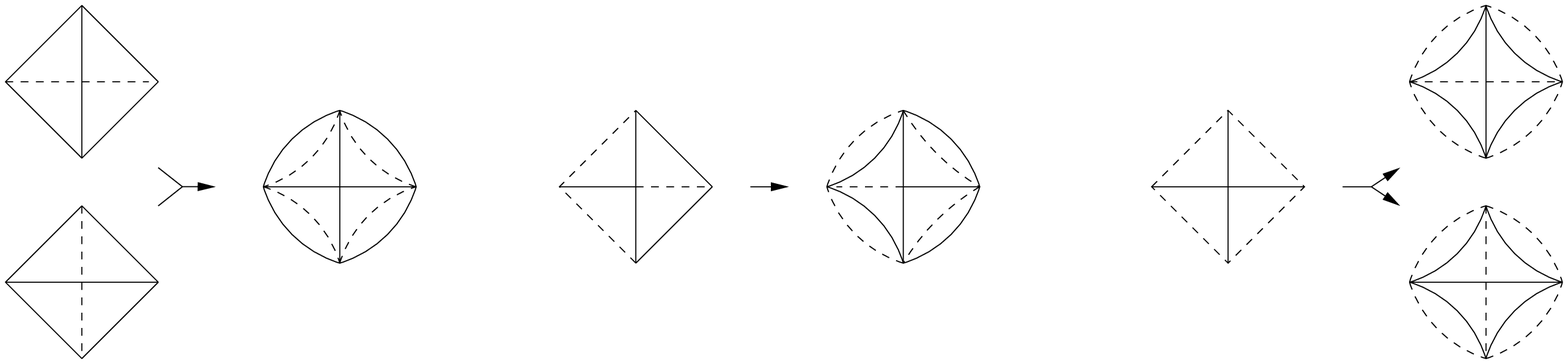,width=10.5cm}\end{matrix}.
\end{equation}

\noindent In the first case of (\ref{bfk}), it turns out that $\overline\mu(M_1)+\overline\mu(M_2)=\overline\mu(M')$, 
in the second case $\overline\mu(M)=\overline\mu(M')$, and in the third case
$\overline\mu(M)=\overline\mu(M'_1)+\overline\mu(M'_2)$. Putting all these together shows (\ref{Spt}).

We do the computation for the first of the three cases of (\ref{bfk}). 
We only take into account the value of $\overline\epsilon$ for the central vertex, which is $-1$.
The local contribution to $\overline\mu(M_1)$ is then 
\begin{equation}\label{asc1}
(ad)^{3/2}(bc)^{1/2}e^{0-1}.
\end{equation}
The local contribution to $\overline\mu(M_2)$\dontshow{asc1} is\dontshow{asc2}
\begin{equation}\label{asc2}
(ad)^{1/2}(bc)^{3/2}e^{0-1},
\end{equation}
and the local contribution to $\overline\mu(M')$ is\dontshow{asc}
\begin{equation}\label{asc}
(abcd)^{1/2}((ad+bc)/e)^{2-1}
\end{equation}
Since $M_1$ and $M_2$ agree
outside of the relevant region, we can group (\ref{asc1}) and (\ref{asc2}) 
together, and say that their joint local contribution is
\[
\big[(ad)^{3/2}(bc)^{1/2}e^{0-1}\big]+\big[(ad)^{1/2}(bc)^{3/2}e^{0-1}\big].
\]
This expression is easily seen to be equal to (\ref{asc}),
 thus showing that $\overline\mu(M_1) + \overline\mu(M_2)=\overline\mu(M')$.
The other two cases of (\ref{bfk}) are similar and can be found in \cite{Spe04}. \qe

Case 3) If the volume between $S$ and $S'$ is a 1/2 octahedron, then we have the following local situations\dontshow{3ca}
\begin{equation}\label{3ca}
\begin{matrix}\epsfig{file=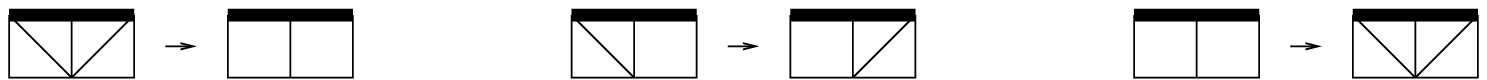,width=10cm}\end{matrix},
\end{equation}
where the thick lines denote the boundaries of $\bW$ and $\bW'$. We do not exclude the case when one or both of 
the points labeled $a$ and $b$ below
lie in $\partial^{(2)}Y$.
As before, we want to show that $Z(\bW,\overline\epsilon,f)=Z(\bW',\overline\epsilon',f')$
where $f$, $f'$ are given respectively by
\begin{alignat*}{1}
\hspace{1.2cm} \begin{matrix}a&e&b\\g&\hspace{.5cm}c\hspace{.5cm}&h\end{matrix} \hspace{1.5cm}&\text{and}\hspace{1.5cm}
\begin{matrix}a&\;\;ab/e\;\;&b\\g&c&h\end{matrix}\;\;,\\
\intertext{
and $\overline\epsilon$, $\overline\epsilon'$ are given respectively by
}
\begin{matrix}\epsilon_1&-1/2&\epsilon_2\\?&-1&?\end{matrix} \hspace{1.5cm}&\text{and}\hspace{1cm}
\begin{matrix}\epsilon_1\!-\!1/2&1/2&\epsilon_2\!-\!1/2\\?&-1&?\end{matrix}\;\;,
\end{alignat*}
where $\epsilon_1$ and $\epsilon_2$ are some numbers in $\{0,1/2,1\}$.
Using Lemma \ref{lemUR} and a trick similar to (\ref{eq.lg}), the three cases (\ref{3ca}) reduce to one:
\begin{equation*}
\psfrag{a}{$a$}
\psfrag{b}{$b$}
\psfrag{c}{$c$}
\psfrag{e}{$e$}\psfrag{a+b-e}{$ab/e$}
Z\left(\begin{matrix}\epsfig{file=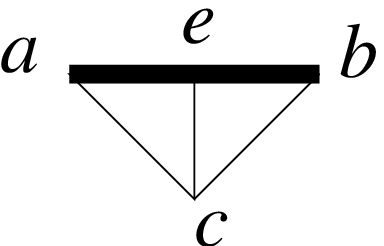,height=.9cm}\end{matrix}\left)=
Z\left(\begin{matrix}\epsfig{file=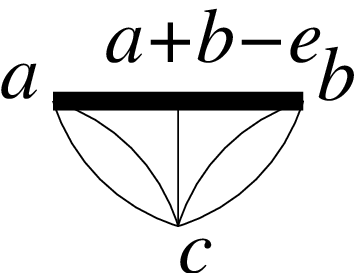,height=1.2cm}\end{matrix}\;\right.\right.\right).
\end{equation*}
This time, the matchings $M$ of the LHS of are in bijection with the matchings $M'$ of the RHS\dontshow{drd}
\begin{equation} \label{drd}
\begin{matrix}\epsfig{file=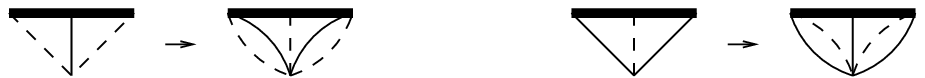,width=7.5cm},\end{matrix}
\end{equation}
and we have $\overline\mu(M)=\overline\mu(M')$. 
We illustrate the computation for the first case of (\ref{drd}). 
The local contribution to $\overline\mu(M)$ is 
$$
a^{-1/2+\epsilon_1}e^{1/2-1/2}b^{-1/2+\epsilon_2}c^{-1/2-1}.
$$
The local contribution to $\overline\mu(M')$ is
$$
a^{0+(\epsilon_1-1/2)}(ab/e)^{-1/2+1/2}b^{0+(\epsilon_2-1/2)}c^{-1/2-1}.
$$
These two expressions are equal, 
therefore $\overline\mu(M)=\overline\mu(M')$. The second case of (\ref{drd}) is similar. \qe
\end{proof}
%Any two sections can be joined by a sequence $S_1, S_2,\ldots $ such that the volume between $S_i$ and $S_{i+1}$ consists of a single 3-cell.
%It follows from Lemma \ref{LsL} that $Z(S)$ is independent of $S$, which 
This last lemma completes the proof of Theorem \ref{tip}.\qed

\section{Applications}

In this section, we prove two surprising facts about the bounded octahedron recurrence which are direct consequences of Theorem \ref{tip}.
Given a subset $S\subset Y$, let $\CH(S):=\semifield^{S\cap \cL}$ denote the set of states of $S$, 
and let $P\CH(S):=\CH(S)/\semifield^\times$ be its quotient by the action of the multiplicative group $\semifield^\times=(\semifield,\cdot)$ of the semifield.
We shall make free use of Remark \ref{Trk}.

Our first application concerns the octahedron recurrence in a $1/4$-octahedron. 
We show that it can be computed using matchings of a hexagon, and deduce from it a hidden 3-fold symmetry.
This is an important example since it has direct connections with the Sch\"utzenberger involution on Young tableaux, and with the theory of $\mathfrak{gl}(n)$-crystals \cite{HK05}.
Matchings of hexagons or ``plane partitions contained in a box'' has been a field of intense activity for its own sake \cite{Pro99},
and it was asked by Jim Propp whether the octahedron recurrence was related to them.

Our second application concerns the octahedron recurrence in a $1/2$-octahedron, 
which we show is equivalent to the recurrence in another domain.

\subsection{The recurrence in a $1/4$-octahedron}

Let $Y=(\R_{\ge 0})^2\times \R$, and consider the domain 
\[
A=\big\{(x,y,t)\in Y\,\big|\,x+y-n\le t\le n-(x+y)\big\}.
\]
Let $S$, $S'$ denote the two equilateral triangles $conv\{(n,0,0),(0,n,0),(0,0,-n)\}$ and $conv\{(n,0,0),(0,n,0),(0,0,n)\}$ respectively.
The bounded octahedron recurrence on $A$ then defines bijections $\phi:\CH(S)\to\CH(S')$ and $\overline\phi:P\CH(S)\to P\CH(S')$.
Let 
\[
\begin{split}
r:S\to S:(n,0,0)\mapsto(0,n,0)\mapsto(0,0,-n)\mapsto(n,0,0)\\
r':S'\to S':(n,0,0)\mapsto(0,0,n)\mapsto(0,n,0)\mapsto(n,0,0)\\
\end{split}
\]
be the rotations of the triangles $S$, $S'$ and let $r^*:P\CH(S)\to P\CH(S)$, $r'^*:P\CH(S')\to P\CH(S')$ be the corresponding pullback maps.
Then we have:
\begin{Proposition}
The maps $\overline \phi$, $r^*$, $r'^*$ defined above satisfy
%Let $\overline\phi:P\CH(S)\to P\CH(S')$ be the map induced by the bounded octahedron recurrence in the domain $A$, 
%induced by the rotations $r$ and $r'$.
%Then one has 
$r'^*\circ\overline\phi=\overline\phi\circ r^*$.
\end{Proposition}

\begin{proof}
Let $f\in \CH(S)$ be a state and $(x_0,y_0,t_0)\in S'$ a point. 
By Theorem \ref{tip}, the value of $\phi(f)$ at $(x_0,y_0,t_0)$ can be computed using matchings of the region 
$W=W(x_0,y_0,t_0)=S\cap\Cmn$ illustrated below:\dontshow{cgo}
{
\psfrag{x}{$\scriptstyle (x_0,y_0,t_0\ts )$}
\psfrag{w}{$W$}
\psfrag{s}{$S$}
\psfrag{sp}{$S'$}
\psfrag{a}{\it \large A}
\psfrag{c}{$\Cmn$}
\begin{equation}\label{cgo}
\begin{matrix}
\epsfig{file=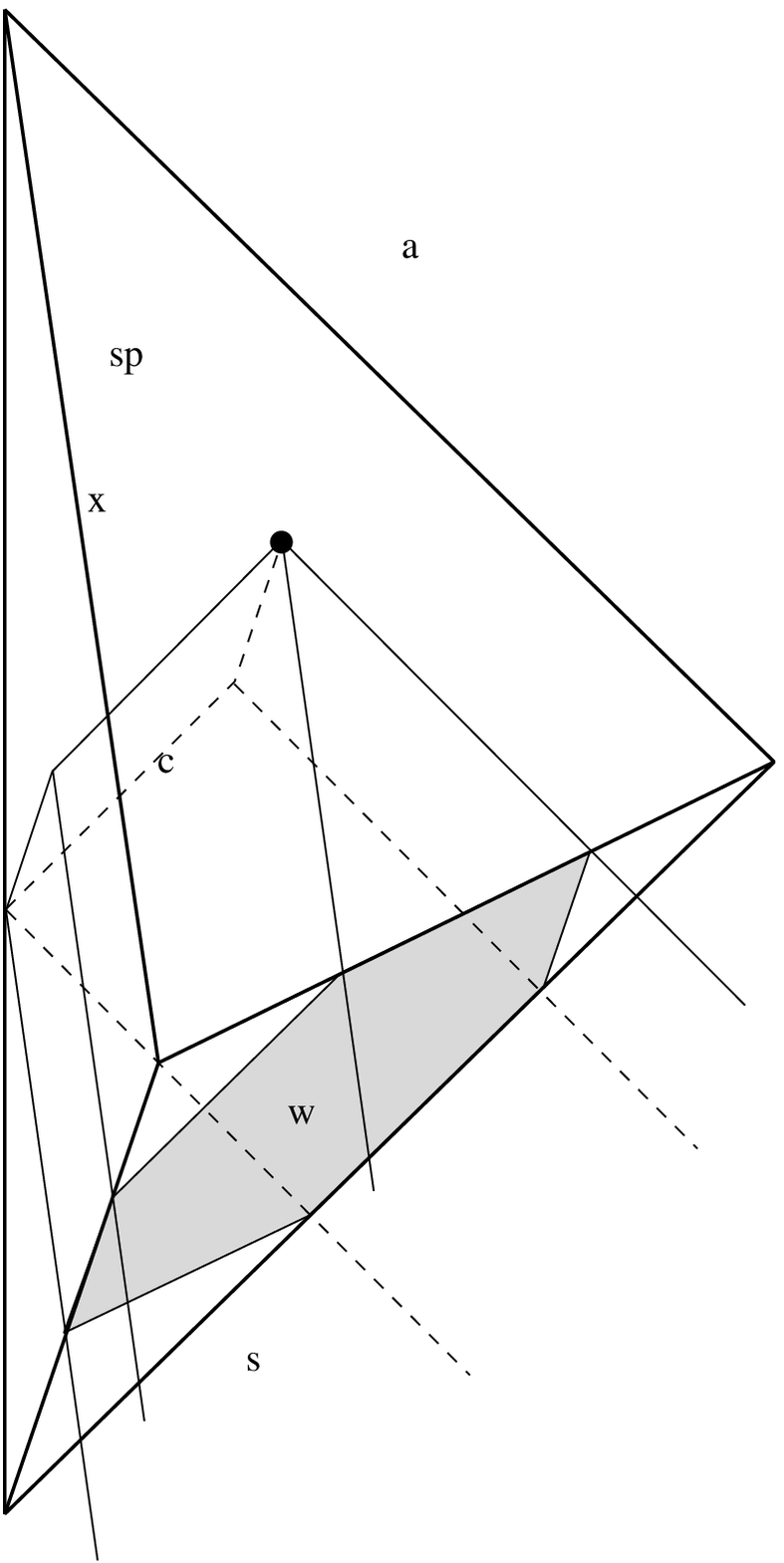,height=5cm}
\end{matrix}
\end{equation}
}

\noindent
A careful analysis of the geometry of (\ref{cgo}) shows that 
\[
W(x_0,y_0,t_0)=\big\{(x,y,t)\in S\,\big|\,x\le n-y_0, y\le n-x_0, -t\le n-t_0\big\}
\]
is a hexagon inscribed in $S$ and that we have $W(r'(x_0,y_0,t_0))=r\big(W(x_0,y_0,t_0)\big)$.
Now we just compute using Theorem \ref{tip}:
\begin{gather*}
\Big(r'^*\circ\phi(f)\Big)(x_0,y_0,t_0)
=
\phi(f)\big(r'(x_0,y_0,t_0)\big)=
c_1\cdot\! \sum_{\substack{M \in \text{matchings}\\ \text{of } W(r'(x_0,y_0,t_0))}}\mu(M,f),\\
\Big(\phi\circ r^*(f)\Big)(x_0,y_0,t_0)
=
c_2\cdot\! \sum_{\substack{M \in \text{matchings}\\ \text{of } W(x_0,y_0,t_0)}}\mu(M,r^*(f))
= 
c_2\cdot\! \sum_{\substack{M \in \text{matchings}\\ \text{of } r(W(x_0,y_0,t_0))}}\mu(M,f).
\end{gather*}
The constants $c_1=f(0,0,-n)^{-1}$ and $c_2=(r^*f)(0,0,-n)^{-1}=f(n,0,0)^{-1}$ are not equal 
but apart from that, the above two expressions are identical.
It follows that $r'^*\circ\overline\phi(f)$ and $\overline\phi\circ r^*(f)$ agree in $P\CH(S')$.
\end{proof}

\subsection{The recurrence in a $1/2$-octahedron}

Let $Y_1=\R\times \R_{\ge 0}\times \R$ and $Y_2=\R_{\ge 0}\times [0,n]\times \R$.
We will compare the bounded octahedron recurrence in the domains
\[
A_1=\big\{(x,y,t)\in Y_1\,\big|\,|x|+|y|+|t|\le n\big\}
\]
and
\[
A_2=\big\{(x,y,t)\in Y_2,\big|\,x+y-n\le t\le y-x+n\big\}.
\]
Let 
\[
\begin{split}
S_1&=conv\{(n,0,0),(0,n,0),(0,0,-n)\}\cup conv\{(-n,0,0),(0,n,0),(0,0,-n)\}\\
S'_1&=conv\{(n,0,0),(0,n,0),(0,0,n)\}\cup conv\{(-n,0,0),(0,n,0),(0,0,n)\}
\end{split}
\] 
be the lower and upper boundaries of the $1/2$-octahedron $A_1$,
and let 
\[
\begin{split}
S_2&=conv\{(0,0,-n),(n,0,0),(0,n,0),(n,n,n)\}\\
S'_2&=conv\{(n,0,0),(0,0,n),(n,n,n),(0,n,2n)\}
\end{split}
\] 
be the corresponding boundaries of $A_2$.
Let $s:S_1\to S_2$, $'s:S'_1\to S'_2$ be the maps
\[
s:\begin{cases}
(n,0,0)&\hspace{-.25cm}\mapsto(n,n,n)\\
(0,n,0)&\hspace{-.25cm}\mapsto(0,n,0)\\
(0,0,-n)&\hspace{-.25cm}\mapsto(n,0,0)\\
(-n,0,0)&\hspace{-.25cm}\mapsto(0,0,-n)\\
\end{cases}
\qquad
s':\begin{cases}
(n,0,0)&\hspace{-.25cm}\mapsto(0,n,2n)\\
(0,n,0)&\hspace{-.25cm}\mapsto(0,0,n)\\
(0,0,n)&\hspace{-.25cm}\mapsto(n,n,n)\\
(-n,0,0)&\hspace{-.25cm}\mapsto(n,0,0)\\
\end{cases}
\]
As before, let 
$\phi_1:\CH(S_1)\to \CH(S'_1)$, $\phi_2:\CH(S_2)\to \CH(S'_2)$ and
$\overline\phi_1:P\CH(S_1)\to P\CH(S'_1)$, $\overline\phi_2:P\CH(S_2)\to P\CH(S'_2)$ 
be the maps induced by the bounded octahedron recurrence, and let
$s^*:P\CH(S_2)\to P\CH(S_1)$, $s'^*:P\CH(S'_2)\to P\CH(S'_1)$ be the pullback maps. 
Then we have
\begin{Proposition}
The above maps satisfy
$s'^*\circ\overline\phi_2=\overline\phi_1\circ s^*$.
\end{Proposition}

\begin{proof}
Let $f\in \CH(S_2)$ be a state and $(x_0,y_0,t_0)\in S'_1$ a point. 
Then by Theorem \ref{tip} we have
\begin{gather*}
\Big(s'^*\circ\phi_2(f)\Big)(x_0,y_0,t_0)
=
\phi_2(f)\big(s'(x_0,y_0,t_0)\big)=
c_1\cdot\! \sum_{\substack{M \in \text{matchings}\\ \text{of } W(s'(x_0,y_0,t_0))}}\mu(M,f),\\
\Big(\phi_1\circ s^*(f)\Big)(x_0,y_0,t_0)
=
c_2\cdot\! \sum_{\substack{M \in \text{matchings}\\ \text{of } W(x_0,y_0,t_0)}}\mu(M,s^*(f))
= 
c_2\cdot\! \sum_{\substack{M \in \text{matchings}\\ \text{of } s(W(x_0,y_0,t_0))}}\mu(M,f),
\end{gather*}
so it's enough to show that $W(s'(x_0,y_0,t_0))=s\big(W(x_0,y_0,t_0)\big)$.
We have two cases, depending on which side of $S'_1$ the point $(x_0,y_0,t_0)$ belongs to.
In both cases, the above identity is checked by 
carefully labeling the vertices in the following figures with their respective coordinates:
{
\psfrag{w}{}\psfrag{y1}{$A_1$}\psfrag{y2}{$A_2$}
\psfrag{c1}{Case 1)}\psfrag{c2}{Case 2)}
\[
\begin{matrix}
\epsfig{file=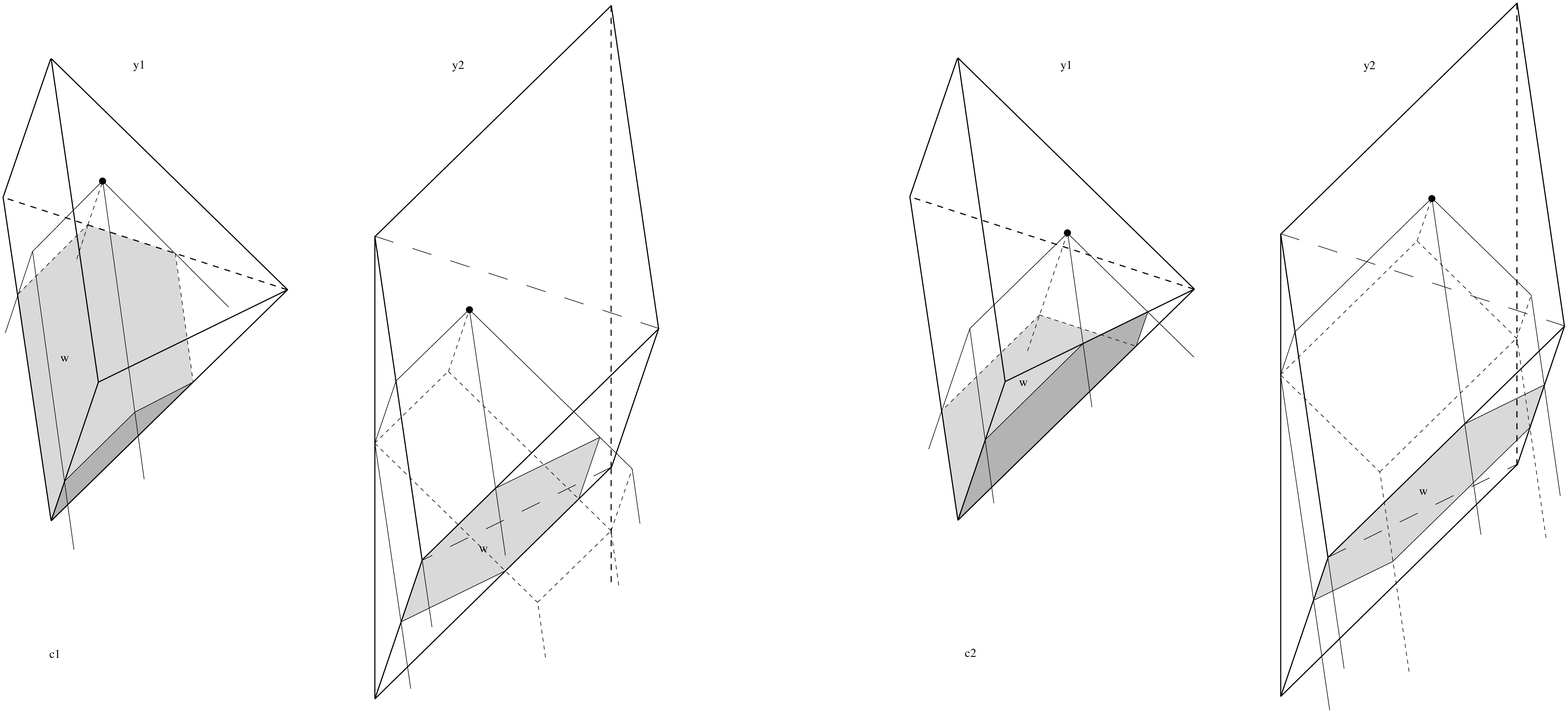,height=5.5cm}
\end{matrix}
\]
}
We leave this task to the diligent reader.
\end{proof}

\section{The cube recurrence}

The cube recurrence, introduced in \cite{Pro01}, is a 3-dimensional recurrence similar to the octahedron recurrence.
It it well known that these two recurrences have a lot of common properties \cite{FZ02, CS04}.
In this section, we introduce a variant which we call the {\em bounded cube recurrence} and make some conjectures
inspired by the bounded octahedron recurrence.

Let $n\ge 1$ be an integer,
let $Y$ be the equilateral prism $\{(x,y,z)\in\R^3\,|\,x\le y\le z\le x+n\}$, and let $\cL:=\Z^3\cap Y$ its set of integral points.
The bounded cube recurrence takes place on the lattice $\cL$ and is given by the formulas
\begin{align*}
\nonumber f(x\ts+\ts 1,y\ts+\ts 1,z\ts+\ts 1&)\ts =\ts
\big(f(x\ts+\ts 1,y,z)f(x,y\ts+\ts 1,z\ts+\ts 1)
               \ts+\ts  f(x,y\ts+\ts 1,z)f(x\ts+\ts 1,y,z\ts+\ts 1)\\ 
               +  &f(x,y,z\ts+\ts 1)f(x\ts+\ts 1,y\ts+\ts 1,z)\big)\big/f(x,y,z)\hspace{.62cm}
\text{if } x\ts <\ts  y\ts <\ts  z\ts <\ts x\ts+\ts n,\\
\nonumber &f(x,y,z\ts+\ts 1)f(x\ts+\ts 1,y\ts+\ts 1,z)/f(x,y,z)          \hspace{.8cm} \text{if } x\ts =\ts y\ts <\ts  z\ts <\ts x\ts+\ts n,\\
\nonumber &f(x\ts+\ts 1,y,z)f(x,y\ts+\ts 1,z\ts+\ts 1)/f(x,y,z)          \hspace{.8cm} \text{if } x\ts <\ts y\ts =\ts z\ts <\ts x\ts+\ts n,\\
\nonumber &f(x,y\ts+\ts 1,z)f(x\ts+\ts 1,y,z\ts+\ts 1)/f(x,y,z)          \hspace{.8cm} \text{if } x\ts <\ts y\ts <\ts z\ts =\ts x\ts+\ts n,\\
\nonumber &f(x,y,z\ts+\ts 1)f(x,y\ts+\ts 1,z\ts+\ts 1)/f(x,y,z)          \hspace{.8cm} \text{if } x\ts =\ts y\ts =\ts z,\\
\nonumber &f(x\ts+\ts 1,y,z)f(x\ts+\ts 1,y,z\ts+\ts 1)/f(x,y,z)          \hspace{.8cm} \text{if } y\ts =\ts z\ts =\ts x\ts+\ts n,\\
\nonumber &f(x,y\ts+\ts 1,z)f(x\ts+\ts 1,y\ts+\ts 1,z)/f(x,y,z)          \hspace{.8cm} \text{if } z\ts =\ts x\ts+\ts n\ts =\ts y\ts+\ts n.
\end{align*}
It is probably best understood by drawing the following pictures
{
\psfrag{a}{$a$}
\psfrag{b}{$b$}
\psfrag{c}{$c$}
\psfrag{d}{$d$}
\psfrag{e}{$e$}
\psfrag{f}{$f$}
\psfrag{g}{$g$}
\psfrag{b+e-g}{$be/g$}
\psfrag{c+f-g}{$cf/g$}
\psfrag{max(a+f,b+e,c+d)-g}{$(af+be+cd)/g$}

\[
\begin{matrix}
\epsfig{file=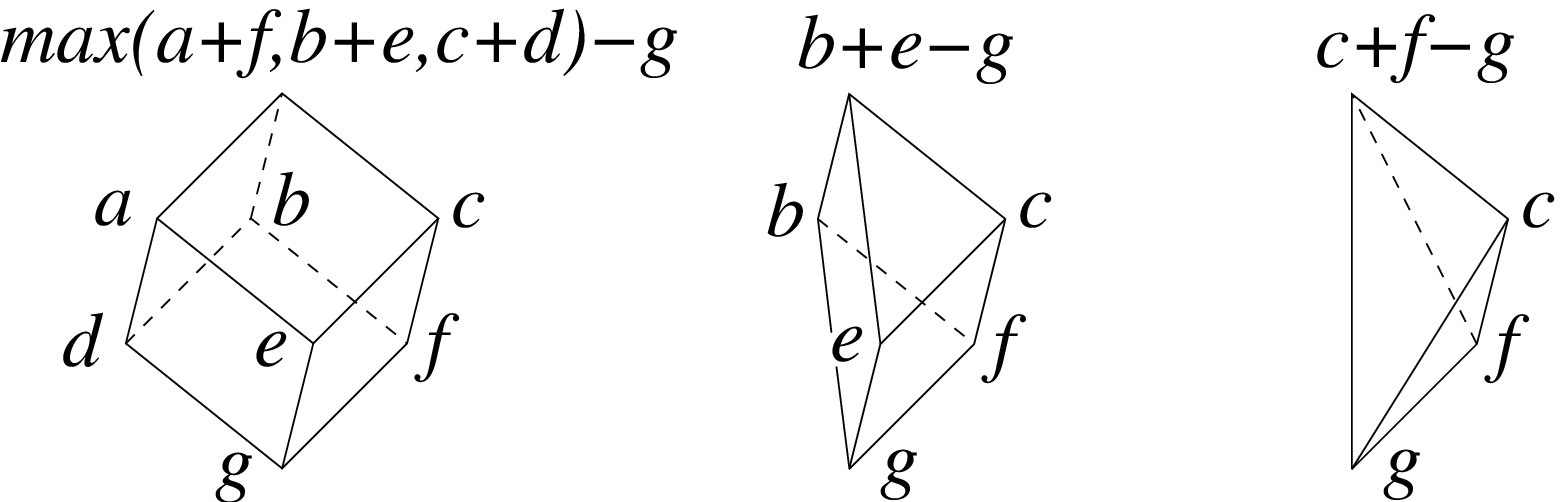,height=2cm}
\end{matrix}
\]
}

\noindent
where the vertical direction represents the vector $(1,1,1)$.

A {\em section} is a two dimensional subcomplex $S\subset Y$ whose projection to the triangle $Y/(1,1,1)\R$ is a homeomorphism.
As in the case of the cube recurrence, we shall say that a point $(x_0,y_0,z_0)$ is in the future of $S$ if there exists $t\ge 0$
such that $(x_0-t,y_0-t,z_0-t)\in S$.
Given a state $f:S\cap \cL\to\semifield$, the bounded cube recurrence then provides an extension of $f$ to all the points of $\cL$ in the future of $S$.

Given such a point, let $\Cmn=\Cmn(x_0,y_0,z_0):=\{
(x,y,z)\in Y|y_0-n\le z\le z_0,z_0-n\le y\le y_0, x_0\le z\le z_0
\}$ be the region illustrated below
{
\psfrag{xyz}{$(x_0,y_0,z_0)$}
\psfrag{y-n,z-n,x}{$(y_0-n,z_0-n,x_0)$}
\psfrag{Y}{\large \it Y}
\psfrag{C}{$\Cmn$}
\psfrag{xyy}{}\psfrag{x,x,z}{}\psfrag{z-n,y,z}{}\psfrag{xxx}{}\psfrag{z-n,z-n,z}{}\psfrag{y-n,y,y}{}\psfrag{y-n,x,x}{}\psfrag{y-n,z-n,y}{}\psfrag{z-n,z-n,x}{}

\[
\begin{matrix}
\epsfig{file=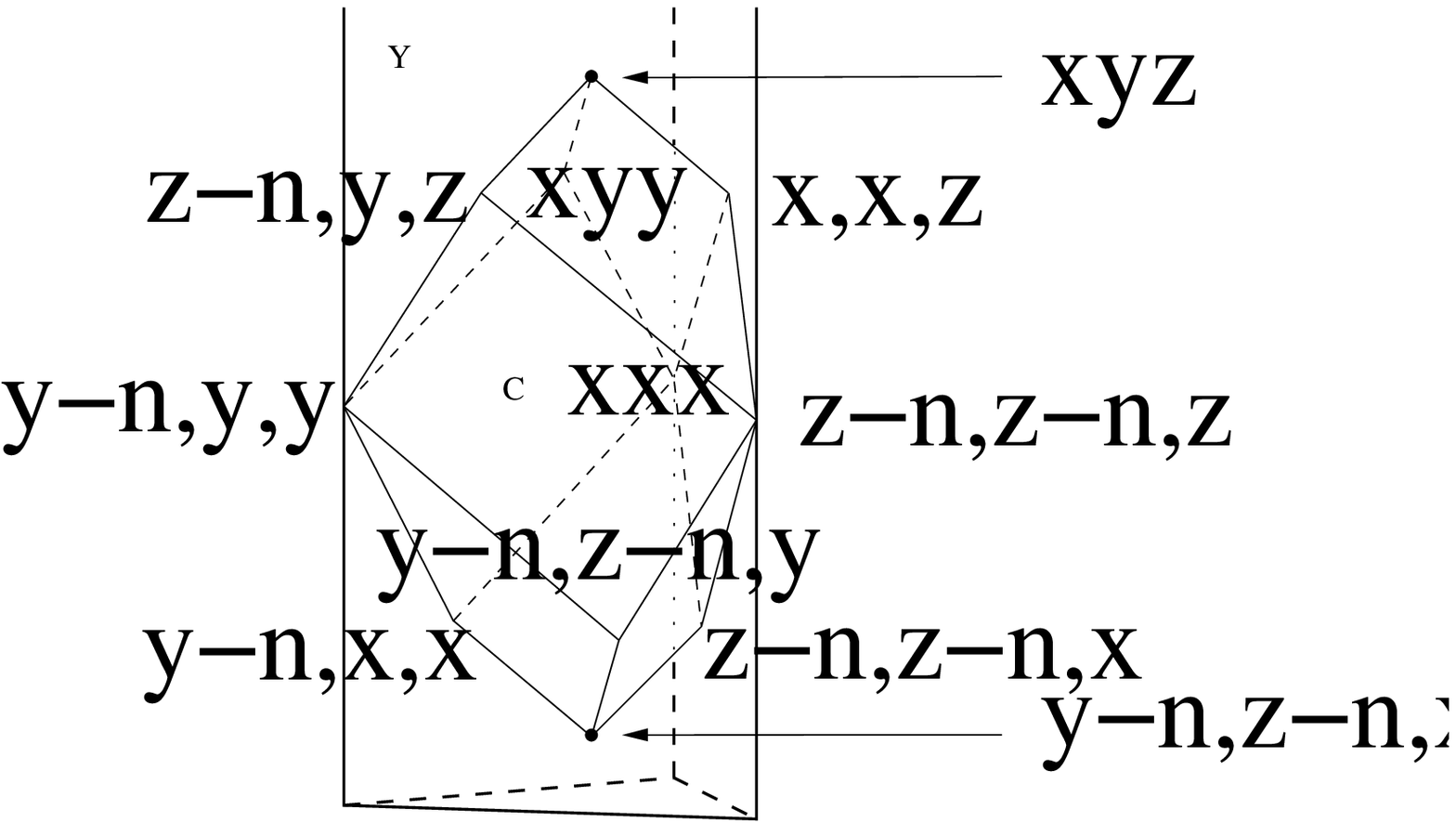,height=4cm}
\end{matrix}\hspace{2cm}
\]
}

\noindent
and let $W$ be the closure of $S\cap \big(\overset{\circ}{\Cmn}\cup\{(x_0,y_0,z_0)\}\cup\{(y_0-n,z_0-n,x_0)\}\big)$.

\begin{Conjecture}\label{cj}
Let $S$ be a section, $f$ a state, $(x_0,y_0,t_0)$ a point in $\overset{\circ}{Y}\cap\tubelattice$ in the future of $S$, 
and assume that the cone $\Cmn$ meets $S$ in at least one point. 
Then there exists a formula of the form
\begin{equation*}
f(x_0,y_0,t_0)=c\cdot\sum_{G \in \text{\rm groves on } W}\mu(G).
\end{equation*}
that computes the output of the bounded cube recurrence is terms of $f|_S$.
We refer the reader to \cite{CS04} for the definition of groves.
\end{Conjecture}

As a corollary of the above conjecture, we have:

\begin{Corollary}
Let $f:\cL\to\semifield$ be a function satisfying the bounded cube recurrence, then there is a constant $c$ such that for any point $(x,y,z)\in \cL$ 
we have the relation
\begin{equation*}
f(x,y,z)=c\cdot f(y-n,z-n,x).
\end{equation*}
\end{Corollary}

%\begin{thebibliography}{E-G-S}
%
%\bibitem{FG05}
%V. V. Fock and A. B. Goncharov,  Moduli spaces of local systems and higher Teichm\"uller theory. arXiv: math.AG/0311149.
%
%\bibitem{us1}
%A. Henriques and J. Kamnitzer, Crystals and coboundary categories. arXiv:math.QA/0406478.
%
%\bibitem{us2}
%A. Henriques and J. Kamnitzer, The octahedron recurrence and gl(n) crystals. axXiv:math.QA/0408114.
%
%\bibitem{KTW04}
%A. Knutson, T. Tao and C. Woodward, A positive proof of the Littlewood-Richardson rule using the octahedron 
%recurrence. arXiv:math.CO/0306274.
%
%\bibitem{Pos}
%O. Gleizer and A. Postnikov, Littlewood-Richardson coefficients via Yang-Baxter equation. \textit{Internat. Math. Res. Notices} (2000),  no. 14, 741--774.
%
%\bibitem{RR86}
%D.P. Robbins and H. Rumsey, Determinants and alternating-sign matrices.  \textit{Advances in Math.} \textbf{62} (1986), 169--184.
%
%\bibitem{Spe04}
%D. Speyer, Perfect matchings and the octahedron recurrence. arXiv:math.CO/0402452.
%\end{thebibliography}
\bibliography{../../main}
\bibliographystyle{plain}

\end{document}